\documentclass[a4paper,12pt]{article}
\usepackage[a4paper, total={7in, 8in}]{geometry}
\usepackage[utf8]{inputenc} 
\RequirePackage[OT1]{fontenc}

\RequirePackage{amssymb,amsmath,amsthm,amsfonts,bbm,mathabx,nicefrac}
\RequirePackage{hyperref,url,booktabs}

\usepackage{natbib}
\usepackage{microtype}      
\usepackage{float}
\usepackage[english]{babel}

\usepackage{graphicx}
\usepackage{comment}
\usepackage{verbatim}
\usepackage{xcolor}
\usepackage{enumitem}
\usepackage{multirow}

\newcommand{\indep}{\perp \!\!\! \perp}

\numberwithin{equation}{section}

\theoremstyle{plain}

\newtheorem{theorem}{Theorem}[section]
\newtheorem{lemma}[theorem]{Lemma}

\newtheorem{algorithm}{Algorithm}[section]
\newtheorem{assumption}{Assumption}
\newtheorem{remark}{Remark}[section]

\theoremstyle{definition}
\newtheorem{definition}{Definition}[section]
\newtheorem{claim}{Claim}[section]

\DeclareMathOperator*{\argmin}{arg\,min}

\title{Model-free Bootstrap and Conformal Prediction in Regression: Conditionality, 
Conjecture Testing, and Pertinent Prediction Intervals}
\author{Yiren Wang\footnote{Department of Mathematics, 
Univ.~of California, San Diego;  email: yiw518@ucsd.edu}
\and Dimitris N. Politis\footnote{Department of Mathematics and Halicio\v{g}lu Data Science Institute, Univ.~of California, San Diego;  email: dpolitis@ucsd.edu}}
\date{}
\begin{document}

\maketitle

\begin{abstract}
	Predictive inference under a general regression setting is gaining more interest in the big-data era. In terms of going beyond point prediction to develop
prediction intervals,  two main threads of development
are  conformal prediction and Model-free prediction. Recently, 
a new conformal prediction approach was  proposed that exploits the same uniformization procedure as in the well-known Model-free Bootstrap.
Hence, it is of   interest to compare and further investigate the performance of the two methods. 
In the paper at hand, we contrast the two approaches via theoretical analysis and numerical experiments with a focus on conditional coverage of prediction intervals. 
We discuss  suitable scenarios for applying each algorithm, underscore the importance  of 
conditional vs.~unconditional coverage, and 
show that, under mild conditions, the Model-free bootstrap yields prediction
intervals with guaranteed better conditional coverage compared to quantile estimation. We also 
extend the concept of `pertinence' of prediction intervals   to the nonparametric regression setting, and give concrete examples where its importance emerges under finite sample scenarios. Finally, we define the new notion of `conjecture testing' that is the analog of
hypothesis testing as applied to the prediction problem; we also devise a 
modified conformal score to allow  conformal prediction
to handle one-sided `conjecture tests', and compare to the Model-free bootstrap.\\


\end{abstract}

\section{Introduction}\label{sec:intro}
Regression is a ubiquitous tool in   statistical literature
and practice. To discuss it, let $X\in \mathbb R^d$ denote the random regressor/covariate, and $Y\in\mathbb R$ the (univariate) response associated with $X$. 
A   general form of a regression setup is the following: 
the relationship between the response $Y$ and covariate $X$ is described via their (assumed continuous) joint cumulative distribution function (CDF)   $F$ which we
denote $(X,Y)\sim F$. In particular, consider the setup where the observations
\begin{equation}\label{eq:regression}
	(X_1,Y_1), \cdots, (X_n,Y_n) \ \mbox{are i.i.d. from} \  F 
\end{equation}
were i.i.d.~is short for 
independent, identically distributed.
This general regression setup has allowed 
researchers to study valid statistical procedures for inference and prediction problems. 

The focus in this paper is predictive inference, specifically, the aim
to generate a valid prediction interval (PI) 
for the response  $Y_{f}$ associated with a future covariate $x_f$ of interest; the PI
should have  a predetermined coverage level/size. 
Formally, we call $\mathcal C_{1-\alpha}(X_{f})$ an `oracle' PI of size $1-\alpha$ at a new covariate $X_{f} = x_f$, if under (\ref{eq:regression}),
\begin{equation}
	P(Y_{f}\in \mathcal C_{1-\alpha}(X_{f})|X_{f} = x_f) = 1-\alpha.
\end{equation}
In practice, the oracle PI is unknown; it  has to be estimated through a statistical procedure,
i.e.,  given data $\{(X_i,Y_i)\}_{i=1}^n$ the statistician produces 
an estimate $ \mathcal{\widehat C}_{1-\alpha}(X_{f})$ such that $P(Y_{f}\in \widehat C_{1-\alpha}(X_{f})|X_{f} = x_f) \approx 1-\alpha$. This estimated PI is called  {\it conditionally valid} by some authors, see e.g., \citet{nonparam_conformal};  we will further discuss the issue of conditionality in prediction intervals in Section   \ref{sec: conditionality}.

Under the general  regression setup of \eqref{eq:regression}, there are at least three 
avenues in constructing prediction intervals, namely: naive quantile estimation, 
conformal prediction, and the Model-free Bootstrap (MFB);
they are all discussed in detail in Section \ref{se:MF reg}.

To give a short preview: quantile estimation (QE) amounts to estimating the conditional distribution 
of $Y_f$ given $X_f=x_f$ and using the associated quantiles   to form
the prediction interval. The term `naive' is due to the fact that using the
quantiles of the estimated conditional distribution as if they were exactly
true sweeps under the carpet   the variability of the 
estimated quantiles; this typically  
results in undercoverage, i.e., the size of the QE PI will be
less than then nominal $ 1-\alpha $ in finite samples.

First proposed in \citet{vovk}, conformal prediction (CP, for short) is an ingenious technique for generating prediction intervals for i.i.d. or just exchangeable data  $ Z_ 1 \ldots ,Z_n$; its basic idea is to gauge how well a new data point  $Z_{n+1}$ ``conforms" to the observed sample  $\{Z_i\}_{i=1}^n$ via a self-chosen ``conformity score". A tutorial on conformal prediction can be found in \citet{conformaltutorial}.
There has been a recent surge in extending the conformal prediction idea  to various statistical setups. Notably, \citet{Lei2013}(see also \citet{leijingcp})  extended this idea to the regression model \eqref{eq:regression}, where the newly proposed conformal prediction algorithm constructs $\widehat{\mathcal C}_{1-\alpha}$ that satisfies the following no-undercoverage guarantee:

\begin{equation}\label{eq:CPcoverage}
	\mathbb P (Y_f \in \widehat{\mathcal C}_{1-\alpha}(X_{f})) \geq 1-\alpha.
\end{equation}
The PI is marginally valid as the probability is evaluated unconditionally, based on the randomness of $\{(X_i, Y_i)\}_{i=1}^n \cup (X_f, Y_f)$; its coverage guarantee relies on the exchangeability of the $n+1$ data pairs $(X_i, Y_i)$ as implied by \eqref{eq:regression}. However, the PI generated by the algorithm  does not guarantee conditional validity  as shown in \cite{nonparam_conformal}  where a new algorithm was   proposed to construct a PI with asymptotic “local validity”. In \cite{1905.03222}, the authors proposed a conformalized quantile regression approach which offers further improvement for heteroscedastic data. 
A recent endeavor to extend regression-based conformal prediction to guarantee conditional validity is the so-called distributional conformal prediction  \cite{distributionalconformal}; see our  Section \ref{se:MF reg} where we explain the procedure in detail. A related approach was proposed in \cite{pmlr-v108-izbicki20a}, 
focusing on the split conformal prediction scheme.  Recently, \cite{zhou-muller} generalized the approach from \cite{distributionalconformal} towards random objects in metric space by using transport distance between distributions as conformity score.
Finally, \cite{2408.16381} approaches the problem of quantifying uncertainty for interval estimation with a combination of conformal prediction and the bootstrap.
\color{black}

Alternatively, since its introduction by B. Efron in \cite{efron1979} for i.i.d.~data, the bootstrap procedure has become a longstanding and widely successful technique for quantifying uncertainty in general statistical methods. The fact that the bootstrap --by generating a large number of artificial samples--
can also capture the variability of estimated
quantities inherent in prediction makes it an ideal candidate for predictive inference. \color{black}
To generalize the bootstrap for predictive inference without strict model-based assumptions,  Politis \cite{Politis2015} proposed the {\it ``Model-free prediction principle"} which amounts to 
using an invertible transformation to map a complex, non-i.i.d. dataset, to a dataset consisting
of i.i.d. variables. The motivation is that prediction is easier --even trivial-- with i.i.d. 
data; one can carry out the prediction in the i.i.d. domain and then map back via
the inverse transformation.  Furthermore, the  {\it ``Model-free bootstrap"}(MFB) amounts to
performing the standard bootstrap  in the i.i.d. domain and then map back, yielding
bootstrap PIs in the complex data domain. 
Notably, this general principle bears resemblence to some techniques in the machine learning community, such as the variational autoencoder and diffusion model, for resampling complex datasets. 
\color{black}

In a general regression setting, \cite{Politis2013} first proposed the use of a `local' 
probability integral transform (PIT), i.e. based on locally estimated conditional distributions,  
to map the response data $Y_1, \ldots, Y_n$ to a sample of 
$n$ i.i.d.   variables with uniform (0,1) distribution. In the meantime, the PIT transform also acts as the cornerstone for conformity score construction in \cite{distributionalconformal} and other literature following their approach. Therefore, the two procedures appear to share certain similarities. However, quantile estimation, the bootstrap and conformal prediction are fundamentally different methods that provide distinct types of guarantees. It is important to understand this difference in order to provide practitioners a guide to use these methods appropriately. Thus, our main focus of the paper is to provide an answer to this question, by analyzing the commonalities/differences between QE, CP and MFB, with an emphasis on conditional coverage.

Our paper is organized as follows. In Section \ref{se:2},  we contrast model-free to model-based regression, and discuss 
several different PI constructions under a nonparametric setting.  In Section \ref{sec: asymptotics}, we carry out 
an asymptotic analysis for the conditional coverage probability of the three methods
and show that ---although no method  has a finite-sample    guarantee in terms of conditional coverage--- all three methods satisfy asymptotic conditional validity; see
our Theorem \ref{thm:asymptotics} in Section  \ref{sec:further_analysis}. 
In addition, while we make the case that it is crucial to construct (and evaluate) the PI
conditionally on $X_{f} = x_f$, we delineate three additional possibilities for 
conditioning in Section \ref{sec: conditionality}; we also show that 
in a large-sample setting these three probability measures are reconciled with 
each other.

In Section  \ref{sec:further_analysis}, we further focus on MFB, and show that beyond holding asympototic validity, it has a favorable property, named the {\it asymptotic pertinence property} of PIs; see Theorem 
\ref{thm:pertinence}. Pertinence    improves finite-sample conditional coverage 
by capturing the variability of estimated quantities in the PI construction, a 
feature that is lacking in both QE and CP. 
Also in Section  \ref{sec:further_analysis}, we are able to show theoretically that
MFB PIs have better finite-sample conditional coverage compared to QE, thus 
showing the practical effects of pertinence; see Theorem 
\ref{thm:performance}. 
Unfortunately, no such finite-sample result is available at the moment
to compare MFB to CP. For this reason, we engage in  
a numerical comparison of all three methods in  Section \ref{sec:numerical};
the simulation experiment is very computer-intensive but it is worth it as it offers 
several insights on the practical implementation and performance of the three methods. 

Furthermore, to extend the applicability of prediction intervals for inference, in Section \ref{sec:conjecture} we introduce the new concept of {\it conjecture testing}  that is the analog of hypothesis testing as applied to the prediction problem; we also devise a 
modified conformal score to allow  conformal prediction
to handle one-sided `conjecture tests', and compare to the Model-free bootstrap.

\color{black}
\section{Model-free vs.~model-based regression}
\label{se:2}

\subsection{Model-based regression}
\label{se:MBR}
The study of interval estimation, i.e., confidence intervals, is almost a century old, with its origins dating back to \cite{Neyman}.
The extension to prediction intervals (PI) was soon to follow; see \cite{Patel_PI},  \cite{geisser1993} and \cite{10.1214/21-STS842} for a review. 
In what follows, we will focus on PIs in a regression setting
with  data   $\{(X_i,Y_i)\}_{i=1}^n$; here, $Y_i$ is the univariate response
associated with regressor   $X_i$  that takes values in  $ \mathbb R^d$.
The goal is to construct a PI for the response $Y_f$ that is not yet observed; 
it will be generated in the future, and will be  
associated with regressor value $X_f=x_f.$

\begin{remark} \label{re:1}  \rm We will  assume that the regressors $X_i$ are fixed, i.e., non-random,
	or that they are random but independent of the errors; in the latter case, 
	all the following discussion
	will be understood in terms of the conditional probability
	given $\{X_i=x_i \ \ \mbox{for} \ \ i=1, \ldots , n.\}$
\end{remark}

\vskip .1in
\textbf{Case I:  Model-based
	linear regression  }
Assume the data are generated by the model 
\begin{equation}\label{eq:model-based linear regression}
	Y_i= \beta^{\prime}X_i + Z_i \ \ \mbox{for} \ \ i=1, \ldots , n
\end{equation}
where $\beta \in \mathbb R^d$ is the (unknown) parameter vector, and
the errors $Z_i$ are i.i.d. $N(0, \sigma^2)$.

Letting  $\hat \beta$ and $\hat \sigma^2$  denote the Least Squares (LS) estimator of $\beta$
and $\sigma^2$, an exact $1-\alpha$ PI for $Y_f$ is
\begin{equation}\label{eq:normalPIestLS}
	\widehat {\mathcal C}_{1-\alpha}(x_f) =
	\left( \widehat\beta^{\prime}x_f + t_{n-d,\alpha/2}\widehat\sigma \sqrt{ 1 + x_f^\prime (X^\prime X)^{-1} x_f  }
	, \widehat\beta^{\prime}x_f  - t_{n-d,\alpha/2}\widehat\sigma \sqrt{1 +  x_f^\prime (X^\prime X)^{-1} x_f  } \  \right) 
\end{equation}
where $X$ is the $n\times d$ matrix having $X_i^{\prime}$ as its $i$th row, and
$t_{n-1,\alpha/2}$  is the lower $\alpha/2$ quantile of the $t$ distribution with $n-1$ degrees of freedom.
\\
\\
\textbf{Case II:  Model-based
	nonparametric regression  }
Assume now the data are generated by the model 
\begin{equation}\label{eq:model-based nonparametric regression}
	Y_i= \mu_ {X_i} + Z_i \ \ \mbox{for} \ \ i=1, \ldots , n
\end{equation}
where $\mu_x$ is an unknown (but assumed smooth) function of $x \in \mathbb R^d$,  and
the errors $Z_i$ are i.i.d. $N(0, \sigma^2)$; 
the caveat of Remark  \ref{re:1} still applies. 

The quantity $\mu_x$ can be estimated by any smoothing method, e.g., kernel smoothing, local
polynomials, etc. resulting in the estimator $\hat{\mu}_x$. There is no exact distribution (like the $t$ distribution) to use here,  so we must resort to asymptotic normality
of $\hat{\mu}_x$. The 
analog of
(\ref{eq:normalPIestLS}) is the  PI
\begin{equation}\label{eq:normalPIest nonpar }
	\widehat {\mathcal C}_{1-\alpha}(x_f) =
	\left( \hat{\mu}_ {x_f} - z_{ \alpha/2}\widehat\sigma  V_{ x_f}
	,  \hat{\mu}_ {x_f} + z_{ \alpha/2}\widehat\sigma  V_{ x_f}  \  \right) 
\end{equation}
that only has {\it approximate/asymptotic}  coverage $1-\alpha$. Here, $\widehat\sigma$ is the natural estimate of $\sigma$, and 
$V_{ x_f}$ is a feature of the smoothing method; see    eq. (3.17) of \cite{Politis2015}
for a concrete example. 
The  choice of the smoothing bandwidth is
as important as it is difficult in practice; it  balances the trade-off between the bias
and variance of $\hat{\mu}_ {x_f}$, and    also affects the 
finite-sample coverage of the PI of eq.~(\ref{eq:normalPIest nonpar }).

\color{black}

Model \eqref{eq:model-based nonparametric regression} can be extended to allow for
heteroscedasticity, i.e., letting
\begin{equation}\label{eq:model-based nonparametric regression hetero}
	Y_i= \mu_ {X_i} + \sigma_ {X_i} \epsilon_i \ \ \mbox{for} \ \ i=1, \ldots , n
\end{equation}
where $\mu_x$ and $\sigma_ x$  are unknown (but assumed smooth) functions of $x \in \mathbb R^d$,  and the errors $\epsilon_i$ are i.i.d. $N(0, 1)$.
An analogous PI to  (\ref{eq:normalPIest nonpar }) ensues, having first estimated
$\sigma_ x$ via nonparametric smoothing.  
Nevertheless, this PI, as well as the previously mentioned PIs, i.e.,
\eqref{eq:normalPIestLS} and  (\ref{eq:normalPIest nonpar }), 
are only valid if/when the regression errors are {\it exactly} normal, an
assumption that is often not justifiable. 

Furthermore, it is of interest to go beyond the assumption of a regression equation 
driven by i.i.d. errors since even the most general/flexible such equation, namely
\eqref{eq:model-based nonparametric regression hetero}, is unnessarily restrictive.
For example, \eqref{eq:model-based nonparametric regression hetero} implies that the
skewness and kurtosis of the response $Y_i$ is constant, i.e.,  does not depend on 
$X_i$. However, as argued by \cite{Politis2013}, eq.~ \eqref{eq:model-based nonparametric regression hetero}  can not
be assumed to hold even for   {\tt cps71},  a  dataset available within the {\tt np}
package of R that has served as a workhorse for nonparametric regression  for decades.

Note that in the random regressor case,  a 
model equation such as \eqref{eq:model-based nonparametric regression hetero},
or  any one of the more restrictive homogeneous and/or parametric versions such as 
\eqref{eq:model-based linear regression} or \eqref{eq:model-based nonparametric regression},
all imply  (\ref{eq:regression}) --under the caveat of Remark  \ref{re:1}-- but not vice versa. 
Therefore, assumption  (\ref{eq:regression})  is less restrictive even than general 
nonparametric regression as defined by eq. \eqref{eq:model-based nonparametric regression hetero}. Working under assumption  (\ref{eq:regression}) 
has been termed {\it Model-free regression} by \cite{Politis2015} in the sense that it is 
devoid of a model equation such as \eqref{eq:model-based linear regression}, \eqref{eq:model-based nonparametric regression} or \eqref{eq:model-based nonparametric regression hetero}; this is the
subject of the following subsection.

\subsection{Model-free regression}
\label{se:MF reg}

We now revert to the ``Model-free'' setup of eq. (\ref{eq:regression}) 
which is re-stated for convenience here:
\begin{equation}\label{eq:regression n}
	(X_1,Y_1), \cdots, (X_n,Y_n) \ \mbox{are i.i.d. from} \  F .
\end{equation}
Note that the above implies a random design; it is also possible to conduct Model-free regression
with a fixed design as discussed in  \cite{Politis2015} but we limit the present discussion to 
setup  (\ref{eq:regression n}).
\\

\noindent\textbf {Naive quantile estimation.} The quantile estimation (QE) approach is straightforward;
an oracle PI would be
\begin{equation}\label{eq:qePI oracle }
	(  F_{y|x_f}^{-1}(\alpha/2),  F_{y|x_f}^{-1}(1-\alpha/2)).
\end{equation}
Since   the conditional quantile 
$F_{y|x}^{-1}(p) = \inf \{y\in\mathbb{R}: F(y|x)\geq p\}$ is unknown, 
we may estimate it by 
the appropriate quantile of the CDF estimator, i.e.,
\begin{equation*}\label{nonparam-quantile}
	\widehat F_{y|x}^{-1}(p) = \inf \{y\in\mathbb{R}: \widehat F(y|x)\geq p\}.
\end{equation*}
There is a wide range of selections when it comes to consistent estimators of 
the  conditional CDF; two popular methods will 
be discussed in Section \ref{se: choosing CDF}.  
Having  settled on a  choice of conditional CDF estimator, the
practitioner could plug it in \eqref{eq:qePI oracle }
resulting in the   PI  
\begin{equation}\label{eq:qePI}
	(\widehat F_{y|x_f}^{-1}(\alpha/2), \widehat F_{y|x_f}^{-1}(1-\alpha/2)).
\end{equation}
The above plug-in approach is termed `naive' 
because   it  does not account for the estimation error of the quantiles; hence, 
the PI (\ref{eq:qePI}) is typically characterized by undercoverage
in finite samples. 

To elaborate, whichever consistent 
estimator (call it $\widehat F_{y|x}$) 
you choose, plugging in 
$\widehat F_{y|x}^{-1}(p) $ instead  of $ F_{y|x}^{-1}(p) $ in the
PI construction (\ref{eq:qePI}) is like pretending that $\widehat F_{y|x}$ has no 
variability; the variability becomes negligible in the limit as $n \to \infty$
due to consistency, but not in finite samples.

\vspace{0.3cm}
\color{black}

\noindent\textbf{Model-free bootstrap.}
Consider an arbitrary point predictor $\widehat Y_f$ for 
the response to be observed at future point $x_f$; mathematically,
$\widehat Y_f$ is a function of the observed data $\{(X_i,Y_i)\}$ and $x_f$, i.e., $\widehat Y_f = g(\{(X_i,Y_i)\}_{i=1}^n,x_f)$; $\widehat Y_f$ can also be regarded as
some functional of the estimated predictive distribution, denoted by $\mathcal J(\widehat F_{y|x_f})$. 
In addition, let $\mathcal J(F_{y|x_f})$ be the 
true point predictor that $\widehat Y_f$ is estimating.  
The prediction error $R_f = Y_f - \widehat Y_f$ is termed the “predictive root”
in bootstrap literature;  
through resampling $R_f$, we can approximate the quantiles of the distribution of $R_f$, which then can be used to calculate the bootstrap PI centered at $\widehat Y_f$.

The original bootstrap of \cite{efron1979} is applicable to i.i.d.~data that
are then sampled  with replacement to produce the bootstrap resample. 
\cite{Politis2015} proposed the {\it ``Model-free prediction principle"} which amounts to 
using an invertible transformation to map a complex, non-i.i.d. dataset, to a dataset consisting
of i.i.d. variables. The motivation is that prediction is easier --even trivial-- with i.i.d. 
data; one can carry out the prediction in the i.i.d. domain and then map back via
the inverse transformation.  Furthermore, the  {\it ``Model-free bootstrap" (MFB)} amounts to
performing the standard bootstrap  in the i.i.d. domain and then map back, yielding
bootstrap PIs in the complex data domain. 

In the model-free regression setting, \cite{Politis2013} suggested the use of the
probability integral transform (PIT) --based on locally estimated conditional distributions-- 
in order to transform the response data $Y_1, \ldots, Y_n$ to a sample of 
$n$ i.i.d. variables. 
To elaborate, define the ranks $U_i = F(Y_i|X_i)$ that are i.i.d. $Unif(0,1)$ due to the PIT,
assuming $F_{y|x}$ is continuous in $y$ for all $x$;
hence, this transformation performs a {\it uniformization} of the responses. Meanwhile, we can recover the data $Y_i$ by the identity $Y_i = F^{-1}_{y|X_i}(U_i)$, i.e.,  by utilizing the inverse PIT. 
Asymptotic   validity of the MFB for bootstrap confidence 
intervals was recently shown in \cite{wang2021}.

\begin{remark} 
	\normalfont
	Note that  by letting $U_i = F(Y_i|X_i=x_i)$ where $x_i$ was the realization of $X_i$, and 
	finding that the distribution of $U_i$ conditionally on $X_i=x_i$ is $Unif(0,1)$, i.e., not
	depending on the value of $x_i$, implies that (a) the unconditional distribution of $U_i$
	is also $Unif(0,1)$, and (b)  that $U_i  $  is independent of $X_i$. Hence, 
	the collection of conditional CDFs $\{F_{y|x}: x\in \mathbb R^d\}$ connects two equivalent probability spaces: one is the space of $\{(X_i, Y_i)\}_{i=1}^n$, the other is the space of $\{(X_i, U_i)\}_{i=1}^n$, where $X_i$ and $U_i$ are two mutually independent  sequences.
	In the latter, the regression/prediction problem has been trivialized, bringing it to the
	i.i.d.~setting of Case I of Section \ref{se:MBR}. 
	
\end{remark} 

Here, we propose the following model-free bootstrap algorithm for predictive inference under the distributional setup in (\ref{eq:regression n}). The original MFB algorithm of Politis  \cite{Politis2015}
treated the design points $X_i$ as fixed during the resampling; the algorithm below is similar albeit the
$X_i$'s are also resampled.

\itshape
\begin{algorithm}
	{\bf MFB  prediction interval at future point $x_f$---random regressor case:}\label{algo:MFB} 
	\begin{enumerate}
		\item Use the given data $\{(X_i,Y_i)\}_{i=1}^n$ to estimate the conditional CDF $\widehat F_n(\cdot|X_i)$ and its inverse. 
		\item Calculate estimated ranks $\widehat U_i = \widehat F_n(Y_i|X_i)$. 
		\item Choose a type of predictor $\widehat Y_f = g(\{(X_i,Y_i)\}_{i=1}^n, x_f)$,
		i.e., choose the function $g$,  that will be the center of the prediction interval. The predictive root is defined by 
		$$R_f = Y_f - \widehat Y_f.$$
		\item Create $X_1^*, \ldots , X_n^*$ by resampling the regressors $\{X_i\}_{i=1}^n$, i.e., 
		by  randomly sampling with replacement  from the set  $\{X_1, \ldots , X_n\}$.
		
		\item  
		Create $U_1^*, \ldots , U_n^*$ by resampling the 
		$\{\widehat U_i\}_{i=1}^n$. Let $Y_i^* = \widehat F^{-1}(U_i^*|X_i^*)$. Let $Y_f^*$ be a sample  from distribution $\widehat F(y|X= x_f)$, and 
		$\widehat Y_f^* = g(\{(X_i^*, Y_i^*)\}_{i=1}^n, x_f)$. The bootstrap predictive root is $$R^*_f = Y_f^* - \widehat Y_f^*.$$
		\item The above step is repeated   $B$  times to  create $B$
		replicates of $R^*_f$. Denote by $q_\alpha ^*$ the $\alpha $--quantile of the empirical distribution of
		the $B$ replicates of $R^*_f$.
		The bootstrap prediction interval is then defined as
		\begin{equation}\label{eq:MFB RR PI}
			(\widehat Y_f + q_{\alpha/2}^*, \widehat Y_f + q_{1 - \alpha/2}^* ).
		\end{equation}
	\end{enumerate}
\end{algorithm}

\noindent \rm
Note that there are several possible choices for the function $g$ and the
resulting point predictor $\widehat Y_f$   that serves as the center of the prediction interval. By far the most popular one  is the conditional mean of the estimated conditional CDF of  $Y_f$ given $x_f$. Another popular choice is the conditional median; 
see Politis  \cite{Politis2015} for details.

\color{black}

\noindent \rm
For comparison, we state below the original algorithm model-free bootstrap algorithm, 
i.e., Algorithm 4.4.1 in \cite{Politis2015}: 

\begin{algorithm}
	{\bf MFB  prediction interval at future point $x_f$---fixed regressor case:}\label{algo:MFB2}

	The algorithm is identical to Algorithm \ref{algo:MFB} except that we replace Step 4 by:
	\begin{enumerate}
		\item [4$^\prime$.] Let  $X_i^* = X_i$.  
	\end{enumerate}
\end{algorithm}

\begin{remark} 
	\normalfont
	\cite{Politis2015}   proposed some different variations to the basic Algorithm \ref{algo:MFB2}. 
	For example, the {\it Limit} model-free bootstrap omits step (2) that calculates estimated ranks. Instead, in step 3  the $U_i^*$s are directly sampled from a $Unif(0,1)$ distribution, thus reducing computational cost and simplifying the analysis for proofs. For predictive inference, the author also recommends using the {\it predictive}  model-free approach, which estimates the $t^{th}$ rank $\widehat U_t^{(-t)} = \widehat F^{(-t)}(Y_t|X_t)$ where  $\widehat F^{(-t)}$  is estimated through  the delete-$t$ dataset $\{(X_i,Y_i)\}_{i=1}^n / (X_t, Y_t).$  \cite{Politis2015} also expanded in details the various choices of the CDF/quantile estimation as well as the predictor, see e.g. Section 2.4.2 and Section 4.2 .\color{black}
\end{remark}

\normalfont
\vspace{0.3cm}
\noindent\textbf{Distributional conformal prediction.}
As mentioned in the Introduction, a recent endeavor to extend   conformal prediction to guarantee conditional validity in regression is  distributional conformal prediction (CP).  This was proposed by  \cite{distributionalconformal} who used the aforementioned uniformization transformation --based on a locally estimated PIT-- to calculate the conformity score as follows. 

For a candidate $y\in \mathbb R$ for which we want to test whether it belongs to the PI, $\widehat F_{n+1}(\cdot|x)$ is an estimator for the conditional distribution function $F(\cdot|x)$, based on $n+1$ data $\{(X_i,Y_i)\}_{i=1}^n \cup (X_{f},y)$. Define the sample conditional ranks $\widehat U_i'(X_f, y) = \widehat F_{n+1}(Y_i|X_i)$ for $i = 1,\cdots, n$
and $\widehat U_{n+1}'(X_f, y) = \widehat F_{n+1}(y|X_{f})$. The conformity scores are defined as $\widehat V_i(y) = |\widehat U_i'(X_f, y) - 1/2|$, and the $p-$value for the null hypothesis that $y$ conforms to the data is
\begin{equation*}
	\widehat p(y) = \frac{1}{n+1}\sum_{i=1}^{n+1} \mathbbm{1}(\widehat V_i (y)\geq \widehat V_{n+1}(y)).
\end{equation*}
The $1-\alpha$ prediction interval is defined as the range of $y$ where $\widehat p(y) > \alpha$, i.e. 
\begin{equation*}
	\widehat{\mathcal C}_{1-\alpha}(X_{f}) = \{y : \widehat p(y) > \alpha\}.
\end{equation*}

With certain assumptions, \cite{distributionalconformal} showed that this procedure will guarantee unconditional validity, i.e., a finite-sample guarantee of unconditional coverage not less than $1-\alpha$. Note that under this unconditional setting, the covariate $X_f$ is treated as random, which we argue in Section \ref{sec: conditionality} that this condition can be unsuitable for evaluating performance under a regression setting.

With conditioning considered, the performance guarantee of this proposed method weakened to asymptotic {\it conditional} validity in the following sense:
\begin{equation*}
	\mathbb P (Y_{f}\in \widehat{\mathcal C}_{1-\alpha}(X_{f})|X_{f} = x_f) \rightarrow 1-\alpha, \text{ as } n\rightarrow \infty.
\end{equation*} 
This is a correctness guarantee of prediction intervals and must hold in order for a procedure to be usable. QE, MFB and CP can be proved to satisfy this property under necessary assumptions. However, it does not further differentiate the relative performance between techniques. In Section \ref{sec:further_analysis}, we attempt to find alternative property to back a predictive procedure. Note that the original DCP paper \cite{distributionalconformal} Theorem 3 adopted a different notation $\mathbb P (Y_{f}\in \widehat{\mathcal C}_{1-\alpha}(X_{f})|X_{f} = x_f) \geq  1-\alpha + \mathcal{O}_p(1)$ to make the procedure to appear to be asymptotic no-undercoverage. But essentially it is equivalent to asymptotic conditional validity due to the proof mostly leveraging asymptotic convergence of the conformity scores.

\color{black}


\vspace{0.3cm}

\subsection{On choosing conditional CDF and quantile estimator}
\label{se: choosing CDF}
The conditional CDF and its quantiles plays an important role in all three above mentioned methods. The two quantities are directly related in that we can construct a valid quantile estimator from a CDF estimator, and vice versa. For example, by inverting the  CDF estimator, we get the quantile estimator; to construct the CDF estimator from estimated quantiles, the following relationship is convenient: 
\begin{equation}\label{eq:quantileCDF}
	F_{y|x}(u) = \int_{0}^1 \mathbbm{1}(F^{-1}_{y|x}(\tau) \leq u ) d \tau.
\end{equation}

As one of the fundamental problem in statistics, 
there exists a plethora of estimators with guarantees of consistency. Considering the volume of different types of distributions, methods that do not rely on a particular parametric model for the conditional distribution are preferred in practice. We now list two such methods which will be used for our comparisons of the three PI methods.

\begin{enumerate}
	\item \textbf{Nonparametric CDF estimation.} This method constructs the estimator via nonparametric kernel smoothing as follows: 
	\begin{equation}\label{nonparam-cdf}
		\widehat F(y|x) = \frac{\frac{1}{n}\sum_{i=1}^n W_h(X_i, x)K(\frac{y-Y_i}{h_0}) }{\widebar W_h(x)},
	\end{equation}
	with $W_h(X_i, x) = \frac{1}{h^d}\prod_{s=1}^d w(\frac{X_{i,s} - x_s}{h})$ and $\widebar W_h(x) = \frac{1}{n}\sum_{i=1}^n W_h(X_i,x)$; $w(\cdot)$ is a univariate, symmetric density function with bounded support, and $K$ is a proper CDF; $h$ and $h_0$ are bandwidths that converge to $0$, whose optimal rates depend on asymptotic analysis.
	Under finite sample setting, we choose $h$ and $h_0$ according to existing empirical rule or through cross validation, see Ch. 6, \cite{li_raccine}.
	Nonparametric kernel estimation works well when the dimension $d$ of the covariates is low; 
	with high-dimensional covariates,  the formulation in Equation (\ref{nonparam-cdf}) may lead to poor performance as a result of the curse of dimensionality.
	Instead of using a product kernel, Ch.4, \cite{Politis2015} demonstrated a simplifying approach
	by exploiting the univariate kernel combined with a certain form of distance function $d_{\theta}(\cdot, \cdot)$ in $\mathbb R_d$, and $W_h(X_i,x) = \frac{1}{h}w\left(d_{\theta}(\cdot, \cdot)/h\right)$. This helps when the dimension of $X_i$ is high, or even with  functional-type covariates whose dimension is infinite.
	
	Equation \eqref{nonparam-cdf} is suitable for estimation when $x$ is not a boundary point in the range of covariates. Near the boundary, 
	local linear 
	estimation is preferable; see \cite{DasPolitisregression}
	for details.
	
	\item \textbf{Regression quantiles.} The regression quantile approach can alleviate the curse of dimensionality when estimating conditional quantiles with high-dimensional covariates, and has attracted attention both in statistics and econometrics. It seeks to estimate $ F_{y|x}^{-1}(\tau)$ by searching in the linear span of the covariates 
	$\{X'\beta, \beta \in \mathbb R^d \}$
	such that the expected quantile loss $\mathbb E \rho_\tau (Y - X'\beta)$ is minimized. Here, $\rho_\tau(r) = r(\tau - \mathbbm 1(r<0))$ is the check function.  The rationale behind using regression quantiles is that the global minimizer $q(x)$ to $\mathbb E \rho_\tau(Y - q(x))$ is exactly the $\tau-$quantile 
	$F^{-1}_{y|x}(\tau)$, and we restrict our search for $q(x)$ in the linear span of $x$.
	$\beta$ is estimated by the minimizer of the sample quantile loss:
	\begin{equation}\label{eq:quantileregression}
		\widehat \beta = \argmin_{\beta} \sum_{i=1}^n \rho_\tau (Y_i - X_i'\beta),
	\end{equation}
	which is often calculated via numerical algorithms. The optimization problem in eq. (\ref{eq:quantileregression}) can be further adjusted to estimate $\widehat \beta$ under high-dimensional sparse scenarios, see Chap.15 of \cite{qr}.
	
\end{enumerate} 

\section{A first look at coverage through asymptotic expansion}\label{sec: asymptotics}

The methods introduced in Section \ref{se:2} share some common features. Besides the fact that all of them use the conditional CDF function as an indispensable tool, the distributional conformal prediction and Model-free bootstrap both utilize the PIT to get estimated ranks $\widehat U_t$ as a middle step towards further inference. However, they are also very different by design that result in disparate performance guarantees. 
To further understand this difference, a natural tool is the asymptotic expansion of the coverage probability. In this section, we perform the expansion with a focus on the following conditional coverage
\begin{equation}
	\mathbb P_{x_f} (Y_{f}\in \mathcal {\widehat C}_{1-\alpha}(X_{f})) : =  \mathbb P (Y_{f}\in \mathcal {\widehat C}_{1-\alpha}(X_{f})|X_{f} = x_f). 
\end{equation}
We demonstrate that this traditional approach cannot offer insights to distinguish the performance of these methods. In other words, they share the same form of expansion with incomparable parameters. This calls for new perspectives to analyze their differences, as further detailed in later sections.

As a first step, we perform analysis for the QE method.
\vspace{0.2cm}

\noindent\textbf{Quantile estimation.}
In the example of QE, the estimated conditional quantile $\widehat F^{-1}_{y|x}(p)$ is used to construct the prediction interval. This yields

\begin{equation}\label{eq:anbn}
	\mathbb P_{x_f}(Y_{f}\in \mathcal {\widehat C}_{1-\alpha}(X_{f})) = \mathbb E \left[ F_{y|x_f}\left(F^{-1}_{y|x_f}(1 - \alpha/2) + B_n \right)\right] - \mathbb E\left[ F_{y|x_f}\left(F^{-1}_{y|x_f}(\alpha/2) +A _n \right)\right]
\end{equation}
where $A_n = \widehat F^{-1}_{y|x_f}(\alpha/2) -  F^{-1}_{y|x_f}(\alpha/2)$ and $B_n = \widehat F^{-1}_{y|x_f}(1-\alpha/2) -  F^{-1}_{y|x_f}(1-\alpha/2)$ are quantile estimation errors.

By  a second order Taylor expansion under twice continuous differentiability of the  conditional CDFs \color{black}, expression(\ref{eq:anbn}) equals to
\begin{equation}	\label{eq:coveragebias}
	\begin{split}
		&\left(1-\alpha\right)
		+ f_{y|x_f}\left(F^{-1}_{y|x_f}(1 - \alpha/2)\right)\mathbb E(B_n)    
		- f_{y|x_f}\left(F^{-1}_{y|x_f}(1 - \alpha/2)\right)\mathbb E(A_n) \\
		&+  \frac{1}{2} f_{y|x_f}'\left(F^{-1}_{y|x_f}(1 - \alpha/2)\right)\mathbb E(B_n^2)
		- \frac{1}{2}f_{y|x_f}'\left(F^{-1}_{y|x_f}(1 - \alpha/2)\right)\mathbb E(A_n^2) 
		+ \text{higher-order terms}.
	\end{split}
\end{equation}
With the decomposition $\mathbb E(A_n^2) = Var(A_n) + \mathbb E^2(A_n)$, Equation(\ref{eq:coveragebias}) can also be written in the form of a bias-variance expansion.

A natural next-step for asymptotic expansion is selecting appropriate estimators to optimize  convergence rate. For example, when using nonparametric kernel estimation in the univariate case, the optimal rate for the bandwidth $h$ with non-negative kernels is of  $n^{-\frac{1}{5}}$ to achieve an optimized rate of $\mathcal O(n^{-\frac{2}{5}})$; while the rate of convergence slows down exponentially with respect to the covariate dimension $d$, see Theorem 6.1, \cite{li_raccine}. 
Alternatively, in the case of quantile regression, we can achieve asymptotic efficiency, i.e., a $\sqrt{n}-$convergence rate for the estimated quantiles at a fixed covariate dimension $d$, {\it provided the implied linear regression model for  the quantiles is correct}; 
when the model is false, however,  the bias can be non-negligible and further leads to incorrect coverage. In summary, selecting the appropriate estimator depends on practitioner's discretion on balancing assumption appropriateness and convergence efficiency based on the data at hand. 

While the analysis above provides understanding on the contributing factors to coverage, it can only be used as ways to prove consistency of respective PIs, yet not able to further distinguish the relative performance across these different methods. \color{black} As it turns out, the other two methods have the same form of expansion, which we state as a claim below.

\begin{claim}\label{claim:asymptotics}
	{\it
		Let $e_{total}$ be the ``total estimation error'' regarding the procedures of QE, CP or MFB. Then the first two factors contributing to the coverage bias are the first and second moments of $e_{total}$, i.e., 
		$$\text{Coverage bias} = \mathcal O(\mathbb E_n e_{total}) + \mathcal O(Var_n(e_{total})).$$
		Therefore, if the procedure guarantees $e_{total}$ diminishes in certain ways (e.g. bounded and converge to $0$ in probability), it leads to asymptotic conditional validity of this method. 
		\color{black}}
\end{claim}
The detailed derivation of asymptotic expansion for the other two methods \color{black} are presented in the Appendix section. Note that for CP, our proof involves analyzing the asymptotic behavior at a finer level while \cite{distributionalconformal}  omitted such details by directly assuming the convergence of conformity scores  at a high level .

Claim \ref{claim:asymptotics} provides a way to prove a certain method is asymptotically valid, however it offers little context on  relative performance between different methods. 
In the claim above, the differences on coverage bias between the three methods lie in the coefficients of the first/second order terms, and the format of $e_{total}$, which is an over-simplified expression describing the total estimation error in the corresponding procedure. Both the coefficients and format of $e_{total}$ can become increasingly difficult to analyze as the method becomes more complex. For example, even though QE and CP are both based on the estimated quantiles, $e_{total}$ still takes a different form for CP due to the special procedure within. Moreover, in MFB the total error also includes the error of predictor estimation, and $e_{total}$ of MFB has extra terms comparing with to QE/CP. 
Furthermore, the signs of the coefficients in the expansion can lead to both under/over coverage. 
In summary, using asympototic expansion only guarantees conditional validity of these methods, and cannot further differentiate their performance. \color{black}
We need to study them from new perspectives. In Section \ref{sec: conditionality} and \ref{sec:further_analysis}, we provide one approach to study their difference, through the lens of proper conditioning and asymptotic pertinence.

\section{The conditionality problem}\label{sec: conditionality}
One exciting result about conformal prediction which inspired many recent explorations in this area is the marginal/unconditional validity property as given in eq. (\ref{eq:CPcoverage}). 
The basis of \eqref{eq:CPcoverage} relies on both $(\mathbf X_n, \mathbf Y_n)$ and $(X_f, Y_f)$ being random and exchangeable, or more relevantly, i.i.d. in the model-free regression 
setup \eqref{eq:regression n}; in addition, the coverage probability should be unconditional with respect to the data. 
Discussions about appropriateness of assumptions are mostly focused on the data itself, such as the i.i.d. assumption for the data $(\mathbf X_n, \mathbf Y_n)$ which is standard in the regression literature; or equi-distribution assumption for the future pair $(X_f, Y_f)$, which led to recent advances in conformal prediction that  proposed solutions when there is a distributional shift for $X_f$, see e.g, \cite{1904.06019}. 

What is equally important but less focused on, is the appropriate 
level of conditionality with which we should evaluate the PI's coverage probability. More specifically, at which level should we integrate randomness of data into the probability evaluation. This so-called conditionality problem has been a crucial concept in statistical inference going beyond prediction intervals. Ever since the proposal of the conditionality principle of \cite{conditionality} and its connection to the likelihood principle (see \cite{conditionality_principle}), statisticians have debated on the correctness of these principles and their implications in statistical inference. 
For a review, see \cite{Robins2000}, which also points out that the issue of conditionality does not involve mathematical rigor; instead,  arguments regarding its virtue are  purely a statistical concern. 

In resonance with conditionality for parameter inference problems, 
new findings are recently discovered regarding the impact of conditioning for prediction problems under a regression setup; see e.g., \cite{regression_random_X} and \cite{NCV}, where new tools are also developed for better predictive inference under conditioning. 
In this section, we   study  the evaluation of coverage probabilities under different levels of conditioning, and derive relationships between them. We connect different conditioning
setups with real world scenarios, and suggest the appropriate  level of conditioning under different circumstances. The conditionality problem also  leads to further analysis of the MFB algorithm in Section \ref{sec:further_analysis}.

As a first step, we introduce three probability measures with different levels of conditioning.  For simplicity, they are denoted as $\mathbb P_1$, $\mathbb P_2$, $\mathbb P_3$ respectively:
\begin{equation}	\label{def:prob_measures}
	\begin{split}
		& \mathbb P_1(\cdot) = \mathbb P_{x_f}(\cdot) = \mathbb P(\cdot|X_f = x_f),\\
		& \mathbb P_2 (\cdot) = \mathbb P(\cdot|X_f = x_f; \mathbf X_n),\\
		& \mathbb P_3(\cdot) = \mathbb P(\cdot|X_f = x_f; (\mathbf X_n,\mathbf Y_n) ).
	\end{split} 
\end{equation}

\begin{remark} \rm Note that all three above probability measures include—at the very least---conditioning on $X_f = x_f$.
	By contrast, much of the literature on CP focuses on completely unconditional coverage without even conditioning on $X_f$, and derive finite-sample guarantees that this unconditional coverage is at the desired level.
	
	We argue that lack of conditioning on $X_f=x_f$ is unacceptable when evaluating a PI’s performance. To see why, imagine a toy example where---due to heteroscedasticity---the variance of $Y_f$ is $1$ when $X_f <0$ but it is $10$ when $X_f >0$; let us assume that $X_f$ is Gaussian with mean zero.
	Constructing a PI conditional on $X_f=x_f$ will yield short PIs when $X_f <0$ and large ones when $X_f >0$.
	Constructing a PI that does not depend on $x_f$ will yield a PI that has the same length throughout. The possible results
	are either (a) that the PI is too wide (and overcovers) when $X_f <0$ and too short (and undercovers) when $X_f >0$---with the coverage being about $1-\alpha$ on average, or (b) that the PI is so wide to have correct coverage when $X_f >0$ but then has pronounced overcoverage for $X_f <0$. The disadvantages are obvious: under scenario (a), the PI will nowhere have correct coverage—it will either overcover or undercover; under scenario (b), the PI may just be too wide to be practically useful.

	Even with appropriate conditioning on the future covariate $X_f$, further conditioning on the regressors($\mathbb P_2$ and $\mathbb P_3$) make prediction more difficult, or even impossible for CP, which was recently pointed out in \cite{2402.04859}. On the contrary, we show that the MFB demonstrates flexibility to work with all the conditional settings above, especially the latter two, due to the bootstrap setup. 
	
\end{remark}

We have the following tower property regarding the hierarchical relationship  between $\mathbb P_k : k \in \{1,2,3\}$.
\begin{lemma}(Tower property)\label{lm:tower}
	Let $A \in \sigma\left(\mathbf X_n, \mathbf Y_n, X_f, Y_f\right)$ be an arbitrary measurable event; here $\sigma\left( \cdot  \right)$
	denotes $\sigma$--algebra.  Then, we have that $\mathbb P_2(A) = \mathbb E_{\mathbf Y_n|\mathbf X_n} \mathbb P_3(A)$ and   
	$\mathbb P_1(A) = \mathbb E_{\mathbf X_n, \mathbf Y_n} \mathbb P_3(A)$.
\end{lemma}

\begin{remark} \rm Bootstrap calculations (including the MFB) are 
	typically conducted conditionally  on the data. The so-called `bootstrap world' is 
	governed by probability measure $\mathbb P_3$, and 
	bootstrap validity is typically proven in probability (or almost surely)   with respect to 
	the `real world' probability $\mathbb P $.
	Note that showing validity (or asymptotic validity) under $\mathbb P_3$
	is stronger, implying validity (or asymptotic validity) under $\mathbb P_2$ 
	or  $\mathbb P_1$.
	To see why, assume $\mathbb P_3 (A) \to 1-\alpha $
	in probability; then the expectation of $\mathbb P_3 (A) $
	will also tend to $1-\alpha$  by the bounded convergence theorem, i.e., 
	$\mathbb P_1 (A) \to 1-\alpha$ and $\mathbb P_2 (A) \to 1-\alpha$. Hence, $\mathbb P_3$ is the strongest probability measure for evaluating conditional coverage. We provide two examples where validity under $\mathbb P_3$ can naturally resolve issues with other types of conditioning.
	
	As a first example, apart from the three conditionality problem discussed above, there is another type
	of conditional coverage considered in literature that is termed Probably Approximately
	Correct (PAC) guarantee, see e.g., \cite{pmlr-v25-vovk12} and \cite{10.1145/3478535}. It requires that with high probability,
	\begin{equation}\label{eq:PAC}
		\mathbb P (Y_f \in \hat C_{1-\alpha}(X_f) |\mathbf X_n, \mathbf Y_n)\geq 1-\alpha
	\end{equation}
	The difference in equation \eqref{eq:PAC} and the three conditional probability measures is how we interpret source of randomness when evaluating coverage.  Nevertheless,  conditional coverage guarantee under $\mathbb P_3(\cdot)$ can still lead to equation \eqref{eq:PAC}. Indeed, 
	$$\mathbb P (Y_f \in \hat C_{1-\alpha}(X_f) |\mathbf X_n, \mathbf Y_n) = \mathbb E_{X_f}\mathbb P_3(Y_f \in \hat C_{1-\alpha}(X_f)).$$
	
	Another example is the feedback covariate shift(FCS) problem where the future covariate's distribution can shift according to the observed data. Clearly, the exchangeable property of $\{ (X_i, Y_i)\}_{i=1}^n$ and $(X_f,Y_f)$ no longer holds, and it requires intricate design to adapt conformal prediction to this scenario, see \cite{2202.03613}. On the other hand, if an algorithm's conditional coverage is proved under $\mathbb P_3$, it will be immune to this problem.
	Hence, working under a $\mathbb P_3$ performance
	measure which is more stringent may be generally preferable. 
\end{remark}  

Remark \ref{re:Big} below 
offers a reconciliation of the different viewpoints for the above probability measures in a big data setting.
Consider the following data scenarios.
\vspace{0.3cm}

\noindent\textbf{A. Scenarios with respect to data.} The data $\{(X_i, Y_i)\}_{i=1}^n$ belong to one of the following conditions:
\begin{enumerate}
	\item {  Both $X_i$ and $Y_i$ are random, with i.i.d. $(X_i, Y_i)\sim F$.} 
	\item { $(X_1,\cdots, X_n)$ are fixed covariates; $\{Y_i\}_{i=1}^n$ are considered random with distribution $Y_i|X_i \sim F_{y|X_i}$. }
	\item {The data $\{(X_i, Y_i)\}_{i=1}^n = (\mathbf X_n, \mathbf Y_n)$ are   observations that were drawn from the joint distribution $F$ but are now treated as fixed.}
\end{enumerate}
The above three scenarios are directly coupled with the above probability measures,
i.e., scenario (A.i) corresponds to 
$\mathbb P_i$,  for any $ i = 1,2,3$. In the machine learning community, $\{(X_i, Y_i)\}_{i=1}^n$ are often referred to as training samples. Thus in scenario (A.1), we are interested in the marginal coverage performance averaged across all training samples; in scenario (A.2), we are under a fixed design in the regression setting; while under (A.3), we are interested in the training-conditional coverage.

\begin{remark}
	\label{re:Big} 
  \rm {\bf [Large sample size vs.~limited sample size]}\color{black} 
	By definition of $\mathbb P_3$ in  eq. (\ref{def:prob_measures}),
	the instance-specific $\mathbb P_3$ is a random probability measure which
	inherits its randomness from   the observations $(\mathbf X_n, \mathbf Y_n)$;
	from   lemma \ref{lm:tower} it follows that 
	$\mathbb P_1$ and $\mathbb P_2$ are expectations of $\mathbb P_3$. Under conditions that guarantee the estimated prediction interval converges to the oracle PI, by dominated convergence theorem $\mathbb P_3 (A)$ will converge to $\mathbb P_1(A) $
	for event $A$ of the form $\{Y_f\in \hat C_{1-\alpha}(X_f)\} \in \sigma\left(\mathbf X_n, \mathbf Y_n, X_f, Y_f\right)$  as the  sample size $n\rightarrow \infty$. Therefore, under a \textbf{large sample} \color{black} regime, $\mathbb P_1$ and $\mathbb P_2$ 
	are not very different from $\mathbb P_3$, and are equally suitable for  coverage evaluation. 
	As setups  (A.1) and (A.2) are coupled with $\mathbb P_1$ and $\mathbb P_2$
	respectively, we see that treating the data as random has validity under the large sample regime. 
	On the other hand, when $n$ is finite, $\mathbb P_3$ can be very different from the other two probability measures;  this phenomenon was recently  put forth in \cite{NCV}.  Contrary to the large sample regime, the finite-$n$ case is a \textbf{limited sample}\color{black}  regime where the data should be regarded as fixed observations, i.e., setup (A.3); this suggests the use of $\mathbb P_3$ as being more appropriate for coverage evaluation with limited data.
\end{remark}

To make more concrete the above  large sample vs.~limited sample dichotomy, we provide three leading examples to illustrate     scenarios (A.1), (A.2) and (A.3)
respectively. 
\begin{enumerate}
	\item \textbf{Randomized experiments}. Commonly known as A/B testing, randomized experiments is a prevalent tool for treatment effect estimation. Given a large pool of data $\{X_i\}_{i=1}^n$,
	it involves random selection and assignment of candidates
	into control/treatment groups. Formally put, $\mathcal S_1, \mathcal S_2 \subset \{1,\cdots, n\}$ are disjoint random index sets with $\max(|\mathcal S_1|, |\mathcal S_2| ) \ll n$, and $\{X_{s_1}\}_{s_1 \in \mathcal S_1}$, $\{X_{s_2}\}_{s_2 \in \mathcal S_2}$ are the covariates of treatment/control groups, with responses $\{Y_{s_1}\}_{s_1\in\mathcal S_1}$, $\{Y_{s_2}\}_{s_1\in\mathcal S_2}$ are collected for these covariates. 
	In this case, the data pair $(X,Y)$ in each group should be regarded as random, and $\mathbb P_1$ should be the correct measure for conditional coverage evaluation.
	
	\item \textbf{Regression for categorical data}. Given continuous covariates $X$ and response $Y$, binning (categorization) is a common practice that transforms $X$ into a handful of  categories in order to eliminate unnecessary randomness in the covariates and reduce complexity of data. 
	While $Y$ can still be considered random, binning is a procedure that interpolates between random design and fixed design. It is therefore reasonable to consider $\mathbb P_2$  { when the binning procedure produces fixed covariates}.
	
	\item \textbf{Few-shot learning and modeling of rare events}. The phrase “few-shot learning” is adapted from the machine learning community which describes  scenarios where training data $(X,Y)$ is extremely limited in quantity; this is naturally a limited data scenario. In the case of rare events, $X$ and $Y$ are observations that are extreme values compared to the population distribution. In this case, the coverage probability conditioning on these observations are clearly distinct from one without conditioning. Therefore, $\mathbb P_3$ should be the correct probability measure for coverage.
\end{enumerate}

\

\begin{remark}
	\rm {\bf [Conditionality under prediction.]}
	The conditionality principle was originally proposed in a parameter estimation context; in one of its many versions it states the following: Let $\theta$ be a parameter of interest, then inference with respect to $\theta$, such as evaluating coverage of a CI, should be conditioned on all ancillary statistics, and also relevant subsets; see \cite{Robins2000} for details. 
	
	Under a prediction context, the different levels of conditioning also represent whether we think the data is ancillary and/or relevant subsets with respect to PI coverage. 
	In this regard, $\mathbb P_2$ and $\mathbb P_3$ can be better measures of coverage comparing to $\mathbb P_1$ under certain circumstances. In the above examples, example 2 is related to viewing the covariates as ancillary statistics; while example 3 is related to viewing the entire observations as relevant subsets.
\end{remark}

\noindent
Together with scenarios (A.1), (A.2) and (A.3), we also introduce a second set of scenarios   about the future covariate $X_f$.
\vspace{0.3cm}

\noindent\textbf{B. Distribution assumption for the future covariate.}  
\begin{enumerate}
	\item  $X_f$ is random with distribution $F_X$.
	
	\item $X_f$ is a fixed future design point in $\mathbb R^d$.
\end{enumerate}
(B.1) assumes that the future covariate is randomly chosen; therefore the marginal coverage  $\mathbb P(Y_f \in \mathcal C(X_f))$ is actually a weighted average of conditional coverages: $\mathbb E \mathbb P_{x_f}(Y_f \in \mathcal C(X_f))$.  Being able to achieve a  $1-\alpha$ marginal coverage does not imply an $1-\alpha$ conditional coverage. 
On the other hand, an $1-\alpha$ conditional coverage will guarantee $1-\alpha$ marginal coverage. To this aspect, validity under (B.2) is stronger than that under (B.1).

In recent years, CP has gained increasing interest  in the machine learning community in that it is able to guarantee no undercoverage for a finite number of data under exchangeability;
this is a consequence of (A.1) and (B.1)---see \cite{leijingcp}.
However, CP does not   guarantee conditional validity.
The CP variant   of  \cite{distributionalconformal} offers improvement to \cite{leijingcp} in that it achieves asympototic conditional validity under (A.1) and (B.2).  Yet {\it asymptotically} this method is on par with the performance of naive QE when considering conditional coverage, as demonstrated in section \ref{sec: asymptotics} and appendix \ref{appendix:A1}.
Thus, CP may be  better suited  for large sample scenarios, where data and future covariate can be regarded random.

It has been known that there is no algorithm that has finite sample no-undercoverage guarantees for conditional validity; see Lemma 1 of \cite{nonparam_conformal}. Instead, correctness of a prediction algorithm is demonstrated via asymptotic conditional validity.
As we will show rigorously in the following section, the MFB approach delivers asymptotic coverage validity under any combination of conditions in scenario A and B. 
Granted, this alone does not justify using MFB over QE and CP as both of them are also asymptotically conditionally valid. 
Therefore, besides asymptotic validity, we will introduce two extra properties of
MFB   prediction intervals, namely {\it pertinence} and {\it improved coverage}
over  the other two methods.

\section{Further analysis of Model-free bootstrap}\label{sec:further_analysis}

Recall the two MFB algorithms from Section \ref{se:MF reg}.
Algorithm \ref{algo:MFB} is referred to as the random regressor scheme,
while Algorithm  \ref{algo:MFB2}  is  the fixed  regressor scheme--which can also be
considered when inference is to be made conditionally on the  design covariates $\mathbf X_n$.
Algorithm \ref{algo:MFB} takes advantage of the distributional setup of the covariates, 
and  conforms to the setup where the study of interest is the coverage probability unconditional on the data $\mathbf X_n, \mathbf Y_n$. 
Since both $\mathbb P_2$ and $\mathbb P_3$ assume that the design covariates $\mathbf X_n$ is fixed, Algorithm \ref{algo:MFB2} is tailored for scenarios (A.2) and (A.3), while Algorithm \ref{algo:MFB} is more suitable for (A.1). Asymptotically, the two algorithms are equivalent.

\subsection{Asymptotic pertinence and better coverage guarantees}
In terms of interval prediction algorithms, while asymptotic validity guarantees correctness,  it alone does not justify good finite sample coverage performance. Section 3.6 of \cite{Politis2015} made a clear illustration of this point and defined the notion of ``pertinence"   of prediction intervals. Formally speaking, a  ``pertinent" prediction interval
should  be able to capture both the variability of future response $Y_f$ correctly, as well as the variance due to  estimation.  Although the latter is usually asymptotically negligible, being able to capture it will improve finite sample performance. 
As far as we know, to capture the variance due to  estimation (without resorting
to simplifying assymptions such as normality) requires some form of bootstrap. 

To explain at a gut level the property of pertinence, let us revert to the simple
linear regression with fixed design of Equation \eqref{eq:model-based linear regression}.
If $\beta$ and $\sigma$ were known, an oracle PI for $Y_f$ would be
$$ C^{oracle}_{1-\alpha}(x_f)= \left( \beta^\prime x_f +z_{ \alpha/2 }\sigma,
\beta^\prime x_f -z_{ \alpha/2 }\sigma \right).$$
Comparing the above to the   classical PI given by eq.~\eqref{eq:normalPIestLS}, we see that
$\beta$ and $\sigma$ are---by necessity---replaced by estimated quantities.
The result is that (i) the $z_{ \alpha/2 }$ quantile must be replaced
by $t_{n-d, \alpha/2 }$, and (ii) the inflation factor 
$x_f^\prime(X^\prime X)^{-1\prime} $ must be inserted under the 
square root sign in Equation \eqref{eq:normalPIestLS} to account for the variance of 
$\hat \beta^\prime x_f $.

If the errors are not Gaussian, using the $t_{n-d, \alpha/2 }$
quantile in item (i) above is incorrect---it is here where some kind 
of bootstrap is needed. 
But there are many ways to employ bootstrap here;  
a naive way would be to consistently estimate the distribution of the errors
($G$, say)
by  $\hat G$, and construct the PI
$$ C^{naive}_{1-\alpha}(x_f)= \left(\hat  \beta^\prime x_f +\hat G( \alpha/2 ),
\hat \beta^\prime x_f +\hat G( 1-\alpha/2 )\right).$$
Although asymptotically valid, the above bootstrap PI does not  
capture the variability of $\hat  \beta^\prime x_f $, and thus does not 
address item (ii) above. 
A ``pertinent'' PI makes an attempt to capture the variability of $\hat  \beta^\prime x_f $, as well as employ appropriate quantiles, i.e., address
both items (i) and (ii) above using a carefully crafted bootstrap procedure.

\cite{Politis2015} formulated the notion of ``pertinent"    prediction intervals
in a model-based regression setting such as \eqref{eq:model-based nonparametric regression hetero}. 
In what follows, we extend  the definition of asymptotic pertinence to the model-free regression setup of (\ref{eq:regression n}).
We adopt the usual notation in the bootstrap literature, where $\mathbb P^*$ denote the probability measure in the bootstrap world, i.e. $\mathbb P^*(\cdot) = \mathbb P(\cdot |\mathbf X_n, \mathbf Y_n)$. Similarly, expectation and convergence in distribution in the bootstrap world are denoted by $\mathbb E^*$ and $D^*$, respectively. 

\begin{definition}\label{def:pertinence} \normalfont
	{\bf Asymptotic pertinence of a model-free prediction interval.} 
	A prediction interval generated by the bootstrap predictive root $R_f^*$ (as defined in Algorithms \ref{algo:MFB} and \ref{algo:MFB2}) is asymptotically {\it pertinent} if
	it satisfies all of  the following:
	\begin{enumerate}
		\item Both $R_f$ and $R_f^*$ can be decomposed into the following 
		representations:\\
		$$R_f = \epsilon_f + e_f,$$
		$$R_f^* = \epsilon_f^* + e_f^*$$
		Here, $\epsilon_f =  Y_f - \mathcal J(F_{y|x_f})$ is the non-degenerate variable concerning
		the distribution of future response $Y_f$, and $e_f$ is model estimation error which converges to $0$ in probability;  $\epsilon_f^*$ and $e_f^*$ are their bootstrap analogues. 
		\item $\sup_x \lvert \mathbb P^*(\epsilon_f^*\leq x) - \mathbb P(\epsilon_f \leq x)\rvert \overset{P}{\rightarrow} 0$ as $n\to \infty$. 
		\item  There is a diverging sequence of positive numbers  $a_n$
		such that $a_ne_f$ and $a_ne_f^*$ converge to nondegenerate distributions, and
		$\sup_x \lvert \mathbb P^*(a_n e_f^*\leq x) - \mathbb P(a_n e_f \leq x)\rvert \overset{P}{\rightarrow} 0$ as $n\to \infty$.
		\item  
		$\epsilon_f$ is independent of  $e_f$   in the real world, i.e., 
		in unconditional probability $\mathbb P_1$; similarly, 
		$\epsilon_f^*$  is independent of   $e_f^*$ in the bootstrap world, i.e.,
		in conditional probability $\mathbb P_3$.
	\end{enumerate}
\end{definition}

In section \ref{sec:proofs}, we give closed forms for the decompositions
of part (1) of Definition \ref{def:pertinence}, and show that the MFB algorithm produces asymptotic pertinent prediction intervals, a unique advantage compared to both QE and CP, which do not satisfy the decomposition in Definition \ref{def:pertinence}.
The effect of pertinence will also be shown in section \ref{sec:numerical}, where we run parallel comparison for experiments with  large vs. small sample size $n$, and show that asymptotic validity and pertinence each plays a more important role in the two separate cases.

In addition, we now state a general claim that will be substantiated in what follows. 
\begin{claim}
	\label{claim:5.1} {\it 
		Under mild conditions, the MFB has 
		better coverage probability compared with QE regardless of the type of CDF/quantile estimator one chooses,
		and regardless of which of the three probability measures $\mathbb P_k$  
		one chooses to evaluate coverage. }
\end{claim}

To elaborate,  we will show that 
there exist an integer $N$ and $ \alpha_0 \in (0,1)$, 
such that for all $ \alpha < \alpha_0$ and $n\geq  N$, we have 
\begin{equation}\label{eq:MFBadvantage}
	\mathbb P_k\left(Y_f \in \widehat{\mathcal C}_{1-\alpha}^{MF}(X_f) \right) > \mathbb P_k\left(Y_f \in \widehat{\mathcal C}_{1-\alpha}^{QE}(X_f)\right)
	\ \  \mbox{for any} \ \  k \in\{1,2,3\}.
\end{equation}
Note that for the case of $k = 3$, Equation \eqref{eq:MFBadvantage} holds in probability.

\begin{remark}\rm
	
	The underlying reason behind the above inequality  is that the 
	MFB prediction interval is pertinent, i.e., attempts to capture
	and incorporate the underlying estimation variability, while the
	QE interval does not. For the same reason, 
	we also expect that the MFB prediction interval has provably  
	better coverage as compared to the CP interval. Such a claim can not be supported based on our current technical results,
	i.e.,  the asymptotic expansions in the Appendix, but it is clearly born out in the finite-sample simulations of Section \ref{sec:numerical}. 
\end{remark}

\subsection{Required technical assumptions}\label{sec:proofs}
The idea of proving the property of asymptotic pertinence   and also Equation (\ref{eq:MFBadvantage}) goes as follows. We first consider the setup of (A.3), which is both the natural real-world scenario as well as the correct setup in the bootstrap world. We show that under (A.3) the decomposition of Definition \ref{def:pertinence} holds for MFB, in contrast with QE and CP.
Next, to prove Equation (\ref{eq:MFBadvantage}), we make mild assumption on the estimator $\widehat F_{y|x_f}$, under which the convolution of $\widehat F_{y|x_f}$ with a (asymptotic negligible) Gaussian kernel will lead to more heavy-tailedness of the distribution. 
In the last step, lemma \ref{lm:tower} is invoked to extend our result under (A.3) to  (A.1) and (A.2).

\begin{assumption}  We make the following general assumptions:
	\vspace{.3cm}
	
	(C1). The sequence of conditional CDF estimators $\widehat F_{y|x}$(indexed by $n$)  is consistent to a continuous CDF $F_{y|x}$ in the $\sup$ norm: $$\sup_u \lvert \widehat F_{y|x}(u) - F_{y|x}(u)  \rvert\overset{P}{\rightarrow} 0.$$

	(C2). The MFB is asymptotically valid for the inference problem with respect to $\widehat Y_f$ in the following sense, as $n\rightarrow \infty$:
	\begin{enumerate}
		\item The predictor $\widehat Y_f = g(\{(X_i, Y_i)\}_{i=1}^n,x_f)$ and $\widehat Y_f^* = g(\{(X_i^*, Y_i^*)\}_{i=1}^n,x_f)$ in the bootstrap world both
		satisfy 
		asymptotic normality with matching limiting variance: $\exists \sigma_{\infty}^2 \in \mathbb R_+, \tau_n \rightarrow\infty$ such that 
		\begin{flalign*}
			&\tau_n  \left(\widehat Y_f - \mathbb E \widehat Y_f\right) \overset{D}{\rightarrow} \mathcal N(0, \sigma_{\infty}^2);\\
			&\tau_n  \left(\widehat Y_f^* - \mathbb E^* \widehat Y_f^*\right) \overset{D^*}{\rightarrow} \mathcal N(0, \sigma_{\infty}^2), in\: probability.
		\end{flalign*}
		Where $\mathbb E \widehat Y_f = \mathcal J(F_{y|x_f})+ o(1/\tau_n)$  and $\mathcal J$ is a functional that maps the CDF function to $\mathbb R$; similarly, $\mathbb E^* \widehat Y_f^* = \mathcal J(\widehat F_{y|x_f})+ o_{p}(1/\tau_n)$.

		\item  $\mathcal J(\widehat F_{y|x_f}) \overset{P}{\rightarrow} \mathcal J(F_{y|x_f})$.
	\end{enumerate}

	(C3).  For large enough $n$, $\widehat F_{y|x_f}$ belongs to the following class of CDFs: 
	$$\mathcal F = \{F: \exists \sigma_0, u_0 \in \mathbb R_+, \forall \sigma < \sigma_0, u>u_0, \widebar F(u) < \widebar{F\star\phi_{\sigma}}(u); \widebar F(-u) > \widebar{F\star\phi_{\sigma}}(-u) \}.$$
	Here, $\widebar F(u) = 1 - F(u)$ is the tail distribution function and $\star$ is the convolution operator; $\phi_\sigma$ is the density function of the normal distribution $\mathcal N (0, \sigma^2)$.
\end{assumption}

\begin{remark} \rm 
	The reason we state (C2) as an assumption, rather than formally prove it under additional assumptions, is that the predictor is of practitioner's own choice to make sure it satisfies nice properties, including consistency and asymptotic normality. 
	Such choice is convenient to find. For example, the predictor $\hat Y_f$ can be chosen to minimize the $L_2$ or $L_1$ loss with respect to the distribution of $Y_f$, which means $\hat Y_f$ is the mean/median of the estimated conditional distribution $\hat F(\cdot|x_f)$.  Central limit theorem and the delta method can then be readily used to derive the asymptotic normality condition. There can certainly be other various types of predictors that satisfy asymptotic normality, and it is unlikely we can cover all of them and derive their asumptotic distribution under various kinds of assumptions. Our focus is not on designing such predictors and proving their property, but rather selecting such predictor with the desired property, so that we can leverage the variability of the predictor to achieve higher-order accuracy when using bootstrap to simulate the distribution of the predictive root.
\end{remark}

The above conditions can be easily achieved in general. 
(C1) and (C2) are classic assumptions regarding validity of statistical procedures.
Plus, we make no assumption on the rate of convergence for either the conditional CDF estimator nor the predictor at which we center prediction intervals.  
(C3) is a less intuitive assumption: it implies that the distribution of $Y_{x_f} + \sigma Z$
has heavier tails than that of $Y_{x_f}$ 
when   $Y_{x_f} \sim \widehat F_{y|x_f}$ and $Z\sim \mathcal N (0,1)$
independent of $Y_{x_f}$. 
However, we will show there is in fact a large class of distributions that falls into this class. We start with the following lemma:

\begin{lemma}\label{lm:convolution}
	Let $f(u)$ be a density function defined on $\mathbb R$ where (i) $\exists u_0>0$ such that $f(u)$ is third-order differentiable and $|f^{(3)}(u)|$ is bounded for  $|u|>u_0$.
	(ii) $f(u)$ is convex for $|u|>u_0$, i.e. $f^{''}(u)>0$. Then the CDF of $f$ satisfies (C3).
\end{lemma}

\begin{remark}
	
	\textbf{(a).}Condition (ii) in lemma \ref{lm:convolution} is satisfied by a large class of distributions.
	Since (C3) is a property of the tail distribution function, it is natural to consider distributions that are in the class of regularly varying functions $\mathbf{RV}(\rho)$ at infinity, where the tail CDF $\widebar F(u)$(as well as $\widebar F(-u)$) satisfies
	
	$$\frac{\widebar F(cu)}{\widebar F(u)} \overset{u\rightarrow\infty}{\longrightarrow} c^{-\rho},\, \rho >0;$$
	so that the tail CDF behaves like a power function $u^{-\rho}$ at infinity (multiplied by a slowly varying function $L(u) $). 
	If  $L$, $L'$ and $L^{''}$ behave “nicely”, such that 
	\begin{enumerate}
		\item $L$ is positive and third-order differentiable near infinity;
		\item  $L$, $L'$ and $L''$ all belong to the following class of functions
		$$\mathcal H := \Big\{H: \, \lim_{|u|\rightarrow\infty} \frac{H'(u)}{H(u)/u} = 0\Big\}$$
	\end{enumerate}
	then $\exists C >0$, $f^{''}(u) \sim C u^{-(\rho +3)}L(u) > 0$. The class $\mathcal H$ is a subset of slowly varying functions, see \cite{Resnick2007}.
	
	\textbf{(b).} Another interpretation related to the convolution problem in (C3) is the following. Given $\widehat F_{y|x_f}$ that is consistent to $F_{y|x_f}$ in the $\sup$ norm, then we have that $\widehat F_{y|x_f}\star\phi_\sigma$ is a $C^\infty$ CDF that approximates $\widehat F_{y|x_f}$ in the limit when $\sigma\rightarrow0$. 
	In addition, we have the following result:
	\begin{equation}
		\sup_x |\widehat F_{y|x_f}\star\phi_\sigma(x) - F_{y|x_f}\star\phi_\sigma(x)| \leq \sup_x |\widehat F_{y|x_f}(x) - F_{y|x_f}(x)|,
	\end{equation}
	where the equality holds only when $|\widehat F_{y|x_f}(x) - F_{y|x_f}(x)|$ equals to the $\sup$ norm at every $x$. Then the implication using the  predictive root for bootstrap is two-fold: 
	First of all, the root distribution is more smooth than $\widehat F_{y|x_f}$ itself, and this can be beneficial to the bootstrap procedure, see \cite{10.2307/2336686}. Secondly,  the convoluted CDF has better approximation performance than the plain CDF $\widehat F_{y|x_f}$, this means bootstrapping the predictive root requires less data than only bootstrapping the future response in order to achieve the same performance. This agrees with our numerical findings in Section \ref{sec:numerical}.
\end{remark}

Another important aspect for validating the MFB algorithm is its performance under high-dimensional regime. While we mainly considered fixing the regressor's dimension $d$ in the paper, the assumptions stated here to prove the desired properties of MFB is independent of such restriction. Rather, the assumptions are stated
directly describing the quality of CDF estimators. Even under high-dimensional scenarios, these assumptions can generally hold--albeit the convergence rate may be slower.
Therefore, the prediction interval (PI) generated from MFB will satisfy asymptotic validity as well as asymptotic pertinence under a high-dimensional setup, provided we choose the correct estimation procedure so the assumptions above hold. 

To provide deeper understanding on MFB's performance under a high-dimensional setup, we provide additional simulation results in Section \ref{sec:multivariate_sim}.

\subsection{Results for MFB pertinence and   Claim \ref{claim:5.1}}\label{sec:main_result}

The main theoretical results for our paper are summarized in the following theorems:

\begin{theorem}\label{thm:asymptotics}({\bf Asymptotic validity})
	Assume (C1) with $x = x_f$; then the prediction intervals $\widehat{\mathcal C}(x_f)$ of QE is asymptotically valid.
	Further assume (C1) for all $x$, then CP is also asymptotically valid.
	\begin{equation}\label{eq:validity}
		\mathbb P_k\left(Y_f \in \widehat{\mathcal C}_{1-\alpha}(X_f) \right) \rightarrow 1-\alpha \ \ \mbox{as} \ \ n\to \infty \ \ \mbox{for} \ \ k\in\{1,2,3\}.
	\end{equation}
	(The convergence is in probability sense for $k = 2, 3$.)
	If we also assume $\widehat{Y}_f^* - \widehat Y_f \overset{P}{\rightarrow} 0$, then Equation (\ref{eq:validity}) holds for the MFB as well.
\end{theorem}

\begin{theorem}\label{thm:pertinence}({\bf Asymptotic Pertinence of MFB})
	Assume (C1) and (C2); then,  the predictive root $R_f$ and $R_f^*$ in the Model-free bootstrap is asymptotically pertinent in the sense of Definition \ref{def:pertinence}.
\end{theorem}

\begin{theorem}\label{thm:performance}({\bf Better performance guarantee})
	Assume (C1), (C2)  and (C3); then the Model-free bootstrap has higher coverage comparing to the quantile estimation method in the sense of Equation (\ref{eq:MFBadvantage}).
\end{theorem}
\noindent
The last theorem is in effect verifying our  Claim \ref{claim:5.1}.
Proofs of the three theorems are found in Appendix \ref{appendix:A3}.

\section{Conjecture testing}\label{sec:conjecture}
\subsection{Prediction vs.~estimation}
Consider the two general classes of statistical inference: 
parameter estimation and prediction. While prediction intervals are the analogs 
of confidence intervals in estimation, there is no clearly defined analog 
of hypothesis testing as applied to prediction. We aim to fill this gap here
by defining the new notion of {\it conjecture} testing.
Table \ref{tab:PEP} illustrates these analogies. 
\begin{table}
	\centering
	\begin{tabular}{ |l|c|c| }
		\hline
		Inference type	& Parameter estimation & Prediction\\
		\hline
		\multirow{1}{*}{Point}
		& Point estimation  &Point prediction \\

		\hline
		\multirow{1}{*}{Interval}
		& Confidence interval & Prediction interval  \\
		
		\hline
		
		\multirow{1}{*}{Test}
		&Hypothesis test & {\bf Conjecture test}\\
		\hline
	\end{tabular}
	
	\caption{Analogies between parameter estimation and prediction.}
	\label{tab:PEP}
\end{table}

Conjecture testing is different from hypothesis testing in the following ways.
First of all, conjecture testing attempts to answer a question regarding a future
response $Y_f$ which is an unobserved random variable; by contrast, (frequentist) 
hypothesis testing aims  to answer a question regarding  an unknown fixed
(nonrandom) parameter  $\theta$.
Secondly,  the rejection region of a hypothesis test is derived from the distribution of the target statistic under the null hypothesis;
by contrast, the probability measure by which a conjecture gets tested  is one  of the three mentioned,
i.e., $\mathbb P_k$ (for some $k$),  which remains unchanged when the null changes. 
In addition, contrary to hypothesis testing which often relies on central limit theorems to prove consistency of a particular test, validity of conjecture testing is directly associated with validity of prediction intervals, which essentially requires consistent estimation of the predictive distribution $F_{y|x_f}$. In this regard, as sample size grows, the acceptance region of hypothesis testing will shrink   towards a single point $\theta$, while the acceptance region of conjecture testing will not degenerate.

\subsection{Conjecture testing: definition and discussion}
Hypothesis testing is a fundamental tool in statistical inference. Consider a
parameter of interest $\theta$ whose value is estimated by a statistic $\widehat \theta_n$.
After defining the null hypothesis $H_0$ and the alternative $H_1$,  a 
frequentist  hypothesis test can be conducted   by figuring out
the threshold of the test statistic that ensures $P_{H_0}(\text{reject}\, H_0)\leq \alpha$.
Alternatively, the test can be conducted by computing the $p$-value
and rejecting $H_0$ only if the $p$-value is less than $\alpha$.

By the duality between hypothesis testing and confidence intervals (CI), we can equivalently construct a $1-\alpha$ CI for $\theta$ and  reject $H_0$ if the CI does not include the
$\theta$ value(s) under $H_0$. This duality applies equally to the standard two-sided CIs
that correspond to a two-sided null hypothesis, as well as one-sided CIs ---where one of the
CI bounds are $\pm \infty$---  that correspond to a one-sided null hypothesis.

While prediction intervals are the analogs 
of confidence intervals in estimation, there is no clearly defined analog 
of hypothesis testing as applied to prediction. We aim to fill this gap here
by defining the new notion of {\it conjecture} testing. It all starts by formulating
two conjectures (i.e., hypotheses) regarding the value of the 
future response  $Y_f$, namely  the null conjecture
$C_{null}:\, Y_f \in S_{null}$ vs.~the alternative conjecture
$C_{alt}: \, Y_f  \notin S_{null}$ where $S_{null}$ is a set of interest.

By analogy with hypothesis testing, the 
simplest choices for $S_{null}$ are:
\begin{enumerate}
	\item  $S_{null}=\{ y_0\}$  
	\item  $S_{null}=[ y_0,\infty )$ 
	\item  $S_{null}=(-\infty, y_0]$ 
\end{enumerate}
where $y_0$ is some given value of interest. Case 1 is a point null (with a two-sided
alternative), while Cases 2 and 3 are one-sided tests. The aim is to try to make a decision between the two  conjectures
$C_{null} $ and $C_{alt}$ while controlling the  probability of false rejection of $C_{null}$
to be not more than $\alpha$. The  probability of false rejection
can be measured with probability $\mathbb P_k$  for some appropriate $k$.

As in the parameter estimation paradigm, there is also a natural  duality between conjecture testing and prediction intervals (PI).
For example, in Case 1 of the point null, we could construct a $1-\alpha$ PI for $Y_f$, 
and reject $C_{null}$ at level $\alpha$  if $y_0 \notin $ PI.
If the PI's $1-\alpha$ level is measured with probability $\mathbb P_k$ (for some $k$), 
the  size of the test will also be according to the same $\mathbb P_k$.




It is quite common in practice to be  interested in the behavior of a future response with respect to the predicted value in a particular direction, i.e., leading to
one-sided conjectures and PIs. For example, when excessive risk is associated with response having higher (say) values  than the projection; in this case, the practitioner
would like to find a region of probable values for $Y_f$ that exclude high-risk
extreme values.
Applications can be found in many fields,  such as climate control, system risk monitoring, Value at Risk (VaR), etc.

The above discussion motivates the one-sided tests of Cases 2 and 3.  Focusing on Case  2 ---the other case being similar---
we may start by constructing
a $1-\alpha$ lower one-sided PI of the type   $(-\infty, c_f]$ for $Y_f$, such that $\mathbb P_k\left(Y_f \notin (-\infty, c_f]\right) \approx \alpha$. 
We can then reject $C_{null}: \ Y_f \notin [ y_0,\infty )$
when  
\begin{equation}\label{eq:rejection one-sided}
	y_0 \notin (-\infty, c_f].
\end{equation}
Note that in such a case, $\mathbb P_k(Y_f\geq y_0) \leq \mathbb P_k(Y_f \geq c_f) = \alpha$, i.e., the probability of false rejection of $C_{null}$ will not exceed $\alpha$ as required.

Statistically speaking, under the setup of \eqref{eq:regression n} with i.i.d. data and  future covariate of interest $x_f$, to construct a one-sided PI we need
to  find $c\in \mathbb R$ such that  
\begin{equation}\label{eq:onesidedbound}
	\mathbb P_k (Y_f - \widehat Y_f \leq c) \approx 1-\alpha
\end{equation}
for either $k \in \{1,2,3\}$ depending on the appropriate data scenario.
The inclusion of $\widehat Y_f$ in \eqref{eq:onesidedbound} rules out using the QE, meanwhile we need to carefully design the conformity score for CP to adapt to this one-sided setup.
On the other hand, the Model-free bootstrap  works readily   since 
we can  estimate $c$ by the $1-\alpha$ quantile of the bootstrap predictive root $R_f^*$, i.e., $\widehat c^* = 
\widehat D_{x_f}^{*-1}(1-\alpha)$.
Hence, the threshold $c_f$ appearing in the rejection region \eqref{eq:rejection one-sided} can be
estimated by  $\widehat c_f=\widehat Y_f   +\widehat c^* = \widehat Y_f +\widehat D_{x_f}^{*-1}(1-\alpha)$.

Constructing a one-sided PI -- and inverting it to obtain a one-sided conjecture test -- is also possible using QE. To do that, we can simply use the left/right quantiles of the estimated conditional CDF; for example, a $(1-\alpha)100\%$ one-sided PI based on QE reads: $(\widehat F_{y|x_f}^{-1}(\alpha),\infty)$ or $(-\infty, \widehat{F}_{y|x_f}^{-1}(1-\alpha))$.

\begin{remark}
	Interestingly, the popular CP methodology is not designed to yield one-sided PIs. One might then think that CP can not handle one-sided conjecture tests. Delving deeper, however, it becomes apparent that this limitation is only due to the two-sided nature of the popular conformity score $\widehat V_i(y) = |\widehat U_i'(X_f, y) - 1/2|$. In order to handle one-sided PIs and one-sided conjecture tests via CP we may instead adopt a one-sided conformity measure, such as 
	$\widehat V_i(y) = -\widehat U_i'(X_f, y)$ under which the produced PI will have the form $(c, \infty)$. We used this choice in our real data example in Section \ref{sec:numerical}.
\end{remark}

In Section \ref{sec:numerical}, we give a practical example of applying one-sided conjecture testing to examine accuracy of VaR prediction for high-frequency stock returns data.

\section{Numerical experiments}\label{sec:numerical}
\subsection{Synthetic data with univariate regressor}
\subsubsection{Experiment setup}
We consider the following model for simulations. 
Let $\{(X_i,Y_i)\}_{i=1}^N$ be i.i.d. samples where 
\begin{equation}\label{eq:numericmodel}
	\begin{split}
		& X\sim Unif(0,1),\\
		& Y = sin(\pi X) +\sigma \sqrt{1 + 2X} \epsilon
	\end{split}
\end{equation}
and $\epsilon\sim T_{5}$. We set the future covariate $x_f = 0.5$ to guarantee data balance on the two sides of $x_f$, such that potential boundary effect will not appear in the simulations. The following procedure is conducted to estimate  the conditional coverage probability(CVP) of size $1-\alpha$ prediction intervals averaged on data, i.e., an estimate for $\mathbb P_1(Y_f \in \mathcal C(x_f))$: for each $k \in \{1,\cdots, K\}$, generate $\{(X_i,Y_i)\}_{i=1}^N$ according to \eqref{eq:numericmodel}; then apply QE, CP and MFB on the data to form three $1-\alpha$ prediction intervals $\widehat{\mathcal C}_k^{(QE)}(x_f)$, $\widehat{\mathcal C}^{(CP)}_k(x_f)$ and $\widehat{\mathcal C}^{(MFB)}_k(x_f)$. Also, after fixing $X_f = x_f$, sample from \eqref{eq:numericmodel} $M$ times to get $\{(Y_f)_i\}_{i=1}^M$, and estimate the $k$-th coverage probability by the fraction of $Y_f$ that are in the prediction interval, that is

\begin{equation}
	\widehat{CVP}(\widehat{\mathcal C}_k) = \frac{\sum_{i=1}^M \mathbbm 1\left((Y_f)_i \in \widehat{\mathcal C}_k\right)}{M}
\end{equation}
The collection of $\widehat{CVP}(\widehat{\mathcal C}_k)$ forms an empirical distribution for the coverage probability under model \eqref{eq:numericmodel}. Then we can estimate $\mathbb P_1(Y_f \in \mathcal C(x_f))$ by the sample average $\frac{1}{K}\sum_{k=1}^K \widehat{CVP}(\widehat{\mathcal C}_k) $.

\subsubsection{Parameter tuning}
The aforementioned parameters are set to be the following: $\sigma = 0.2$, $\alpha = 0.05$, $K = 200$ and $M = 3000$. Both the nonparametric  and quantile regression based CDF estimators are considered in our simulations. For the nonparametric estimator, the bandwidth parameters $h$ and $h_0$ need to be selected. We apply the following criteria for selection: within a grid of candidate bandwidths,
the ones with which the uniformized series $\{U_i\}_{i=1}^N$ has the largest $p-$value in terms of the Kolmogorov-Smirnov test for uniformity are selected. As for the quantile regression based CDF, we apply the Barrodale and Roberts algorithm to estimate $F^{-1}_{y|x}(\tau)$ for $\tau$ in the grid $\{0.01i, \, 0\leq i\leq100\}$ and then calculate $F_{y|x}(u)$ using a discretized version of \eqref{eq:quantileCDF}.
In the bootstrap algorithm, we set $B = 1000$.
\subsubsection{Asymptotic performance}
\begin{figure}
	\centering
	\includegraphics[scale = 0.3]{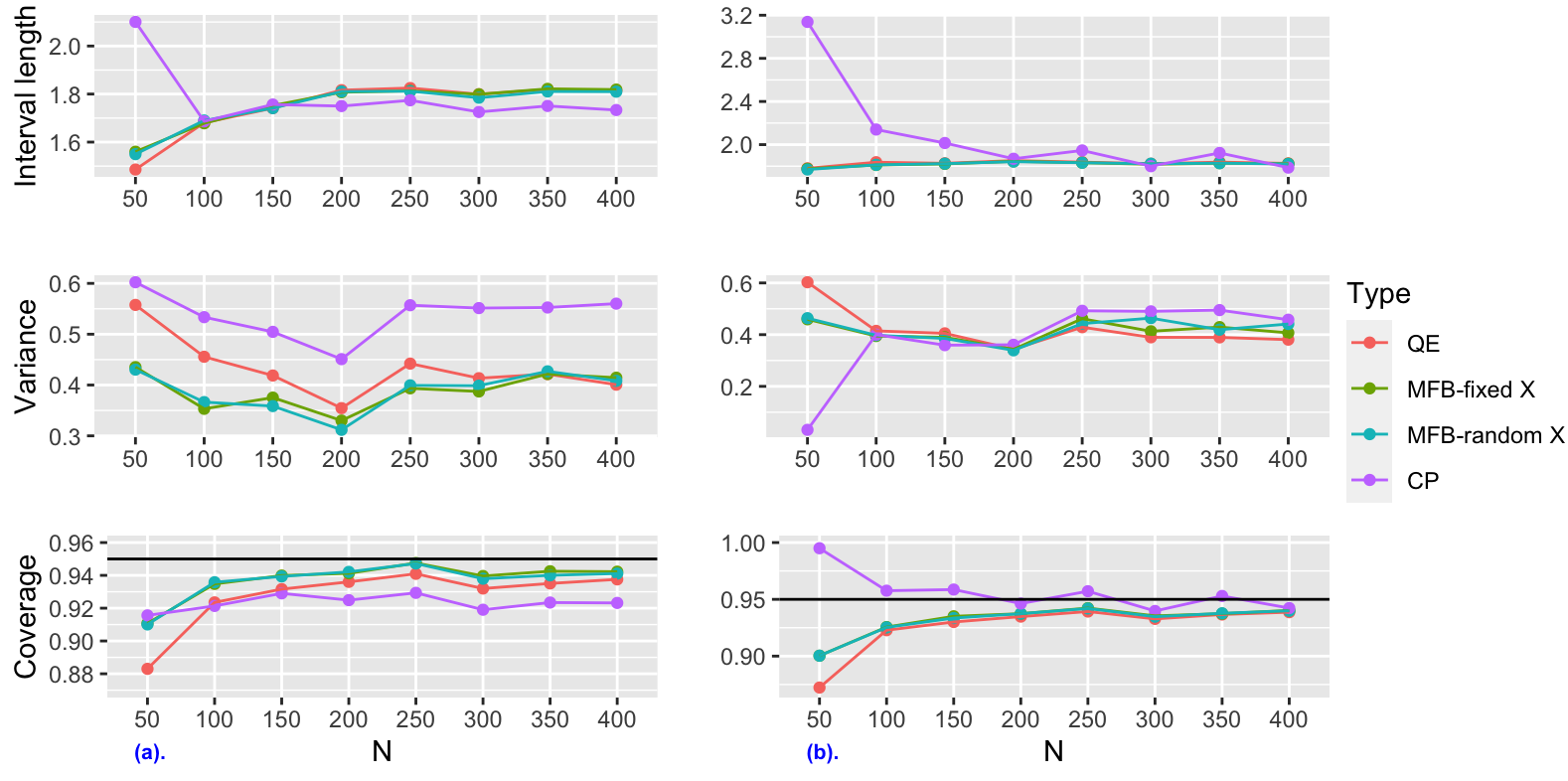}
	\caption{Asymptotic results for prediction intervals. (a): results based on the nonparametric CDF estimator; (b): results based on the quantile regression based CDF estimator.}
	\label{fig:asymptotics}
\end{figure}

We plot the estimated  coverage probabilities, variances of the estimation $Var(\widehat{CVP}(\widehat{\mathcal C}_k) )$(suitably scaled), as well as the average lengths of prediction intervals for different values of $N\in\{50i, i = 1,\cdots,8\}$ in Figure \ref{fig:asymptotics}. 

We analyze the results in the following aspects. 

\textbf{1. Asymptotic validity.} For the nonparametric CDF estimator, QE and MFB exhibit asymptotic validity as the coverage probabilities converge to $95\%$ when $N$ increases; while convergence for CP seems to stabilize around $92\%$. 
For the quantile regression based CDF estimator, all three algorithms converge to the nominal level.

\textbf{2. Coverage comparison.} 
The prediction intervals of MFB have higher coverage than QE for both CDF estimators, especially when the sample size $N$ is small. This agrees with our result in Theorem \ref{thm:pertinence}, and also shows that MFB is able to alleviate undercoverage issue under limited data scenario.
The behavior of CP varies based on the choice of CDF estimator: when the nonparametric CDF is used, CP performance is not as good as the other methods; while for the quantile regression CDF, CP in general has good performance except for the case $N=50$, which we will discuss in further details in the next section.

\textbf{3. Quality of estimation.} The quality of prediction intervals is assessed by the variances of coverage probabilities and also the lengths of the prediction intervals. Differences between the algorithms are most observable from the results that use the nonparametric CDF estimator. 
In terms of the variances, the MFB is consistently lower than QE and CP, while CP has the largest variance among the three methods. This is related to the pertinence property of the MFB that helps tracking the nominal level more consistently than the other methods. 
As for the lengths of prediction intervals, the MFB has equal length compared with QE while delivering higher coverage. On the contrary, CP has an anomaly at $N=50$ where the length is significantly larger; as $N$ increases, the lengths drop below those of QE and MFB, as do the coverage probabilities. 

\subsubsection{Coverage behaviors under limited data scenario}
\begin{figure}[h!]
	\centering
	\includegraphics[scale = 0.32]{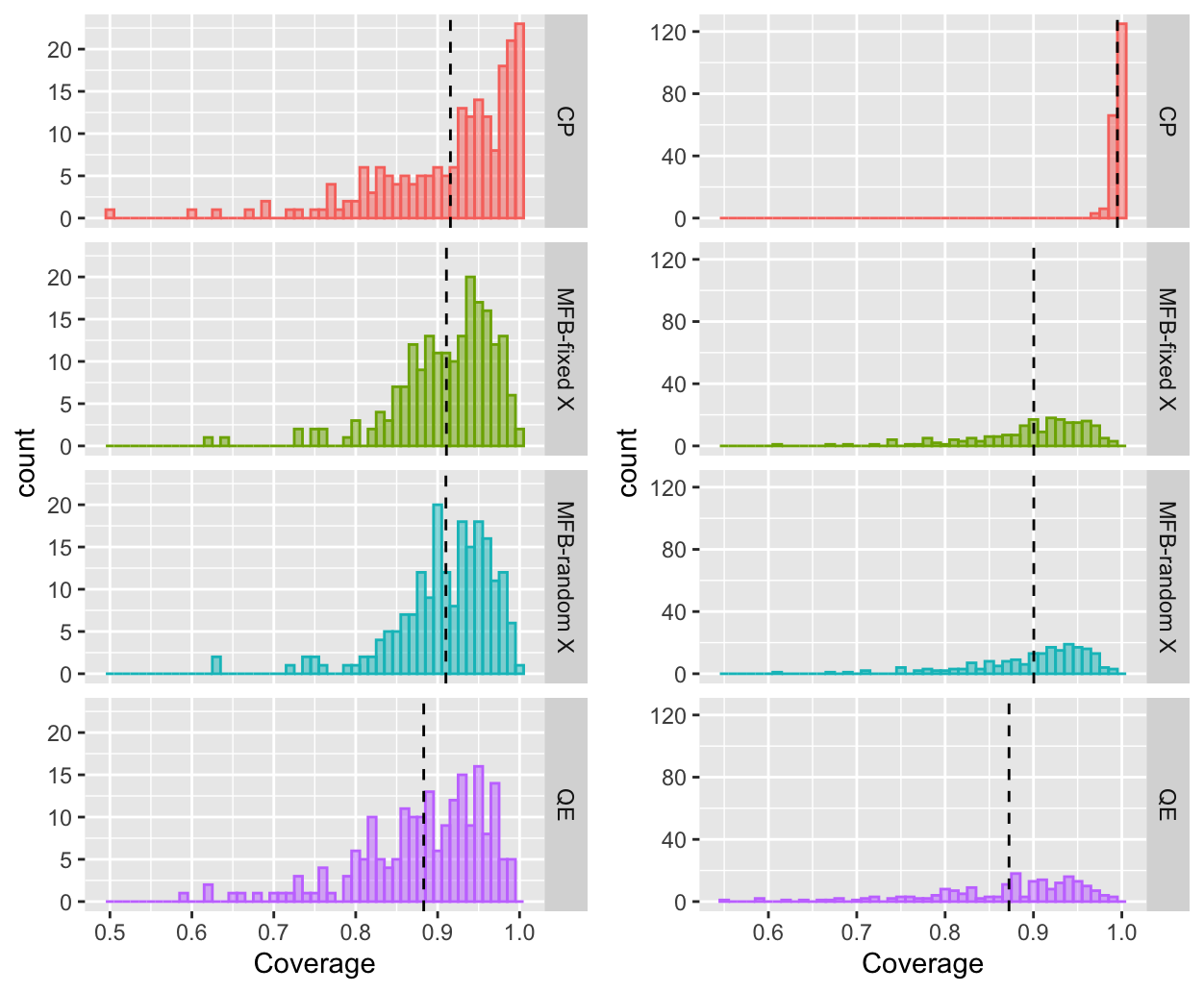}
	\caption{Histogram of coverage probabilities at $N = 50$. Left: results based on the nonparametric CDF estimator; right: results based on the quantile regression based CDF estimator. }
	\label{fig2:histogram}
\end{figure}

The behavior of CP under the limited data scenario ($N = 50$) is abnormal: for both CDF estimators, CP's lengths of  prediction intervals are much higher than competing methods; for the quantile regression based estimator, the coverage is even close to $100\%$. Hence here, we further investigate the results for $N = 50$. In Figure \ref{fig2:histogram} we plot the histograms of the coverage probabilities.
As we see from the graphs, CP will “overcorrect” the prediction intervals to reach the nominal level $95\%$. In the case of using the nonparametric CDF estimator, the average coverage of CP(illustrated by vertical dashed line) is almost identical to that of MFB, but many of the prediction intervals are corrected to reach a $100\%$ coverage, which also cause the interval lengths of CP to be much higher than those of MFB . In the case of using the quantile regression based CDF, this behavior is even more severe, causing the averaged coverage to be almost $100\%$.  
Therefore, the reliability of CP to generate conditionally valid prediction intervals  is questionable under limited data setting. 

\subsection{Synthetic data with multivariate regressors}\label{sec:multivariate_sim}

We also provide an example of applying QE, CP and MFB for a regression model with multivariate regressors. Consider the following model for $\{(\mathbf X_i, Y_i)\}_{i=1}^n$: 
\begin{equation}\label{eq:multivariate linear regression}
	\begin{split}
		&  \mathbf X_{i} \overset{i.i.d.}{\sim} \mathcal N(\mathbf 0, \mathbf I_d) \in\mathbb R^d\\
		& Y_i = \mathbf X_i\boldsymbol\beta + \epsilon_i, \epsilon_i\overset{i.i.d.}{\sim} \mathcal{N}(0,\sigma^2)\\
	\end{split}
\end{equation}
where $\{\beta_i\}_{i=1}^d$ are first one-off generated from $\mathcal N(0,1)$, and then re-scaled so that $||\boldsymbol{\beta}||_2 = 1$. Regarding the CDF and quantile estimators, both non-parametric estimator and quantile regression estimator in section \ref{se: choosing CDF} are considered. In particular, the non-parametric kernel estimator is the multivariate variant of \eqref{nonparam-cdf}, where the distance function is the $L_2$ distance, and the bandwidth $h$ and $h_0$ are chosen such that the uniformized sequence maximizes the p-value of the Kolmogorov-Smirnov uniform distribution test of the sequence. 
We also consider both $L_1$ and $L_2$ predictor in the MFB algorithm. 
For conciseness, we study the case where the sample size $n$ is on the same scale as the covariate dimension $d$, by fixing $d = 50$ and letting $n \in  \{100,150,250, 500\}$. The results of empirical coverage for 90\% level prediction intervals are summarized in Table \ref{tab:simulation_multivariate}.

\begin{table}[h]
	\centering
	
	\begin{tabular}{c|cccc|ccccc} 
		\toprule
		{sample size} & {$QE_{np}$} & {$MFB_{L_1, np}$} & {$MFB_{L_2,np}$} & {$CP_{np}$} & {$QE_{qr}$} & {$MFB_{L_1,qr}$} & {$MFB_{L_2,qr}$} & {$CP_{qr}$} & {$$} \\ \midrule
		{$n=100$} & 91.3\% &90.5\%  & \textbf{90.2\%}  &92.0\%  & 33.7\% &  {\bf 75.8\%}& {\bf 75.8\%} & 100.0\%  \\
		{$n=150$}  &  90.4\%& 90.1\% & \textbf{90.0\%} & 89.0\% & 52.3\% & {\bf 80.9\%} & 80.3\% & 100.0\%  \\
		{$n=250$}  & 87.6\% & \textbf{88.9\%} & 88.8\% & 88.2\% & 65.9\% & {\bf 79.1\%} & 78.9\% & 100\%\\ 
		{$n=500$}  & 89.0\% & {88.3\%} & 88.0\% & \textbf{90.4\%} & 77.8\% & {\bf 80.6\%} & {80.4\%} & 98.4\%\\
		\bottomrule
	\end{tabular}

	\caption{Empirical coverage probability with multivariate regressor from model \eqref{eq:multivariate linear regression}. QE, MFB, CP stand for quantile estimation, model-free bootstrap, and conformal prediction respectively. `np' stands for non-parametric CDF estimator; `qr' stands for quantile regression estimator. } 
	\label{tab:simulation_multivariate}
	
\end{table}

We can derive some useful conclusions from the  results of Table \ref{tab:simulation_multivariate}. First of all, under the high-dimensional setup, the non-parametric CDF estimator has better accuracy for quantile estimation than the quantile regression estimator; this aligns with the result in \cite{NEURIPS2021_9854d7af} stating
that   quantile regression estimator without adjustment can lead to undercoverage for high-dimensional regression models. 
Surprisingly,  results based on the nonparametric CDF estimator are all close to the nominal level under this high-dimensional setup even when the sample size is low. On the other hand, these results do no seem to improve when increasing sample size from $n=100$ to $n=500$; this is due to the curse of dimensionality:  a sample size of $n=500$
is not sufficient  for the asymptotics of CDF estimation to kick in 
under a high-dimensional scenario.

Finally, recall that Claim \ref{claim:5.1} shows that for all $n>$ some $N$ (which may problem specific), the MFB algorithm has higher coverage probability than QE. 
Interestingly, the results of 
Table \ref{tab:simulation_multivariate} do not confirm the claim; this just goes to show
that the $N$ associated with our high-dimensional simulation is quite large, 
certainly larger than 500. Nevertheless,  for both the non-parametric estimator and the quantile regression estimator, the MFB algorithm     leads to better coverage accuracy compared to quantile estimation;
this is   consistently observed regardless as to whether the   quantile estimator exhibits under-coverage or over-coverage. 



\subsection{VaR prediction for intraday stock returns data}\label{sec:sim_2}
We demonstrate an application of conjecture testing on high-frequency intraday stock returns data.  To fix notations, let $S_t$ denote the price process of a stock; $X_t=\log \frac{S_t}{S_{t-1}}$ denotes the log-return $\log \frac{S_t}{S_{t-1}}$ at time $t$. The concept of Value-at-Risk (VaR) is a risk measure in the following sense: for a finite time period $m$ and a pre-specified level $\alpha$, the VaR is the threshold $c$ such that the $m$-period return $T_m = \sum_{t = 1}^{m} X_t$ is guaranteed to surpass $c$ with probability $1-\alpha$. In other words,
$VaR(\alpha) = F_m^{-1}(\alpha)$,
where $F_m$ is the CDF of $T_m$.

This concept can be extended in various ways. 
First of all, we can extend the statistic of interest  beyond the case of total return $\sum_{k = t+1}^{t+m} X_k$, such as the form $T_m = f(X_1,\cdots, X_m)$ where $f$ is a continuous function. Examples include the worst $m-$period cumulative return, $T_m  = \min_{1\leq t\leq m} R_t$, where $R_t = \sum_{k=1}^t X_k$, and the worst one-period return $T_m = \min_{k=1}^m X_k$. These risk-type statistics will induce different forms of $F_m$. These examples are considered in \cite{Bertail2004}. 

Another extension is to find appropriate covariate that can have predictive power for future VaR. For this purpose, we introduce additional dependence on the time parameter $t$ for $T_m$: $T_{t+1,m} = f(X_{t+1},\cdots,X_{t+m})$. Also, let $V_{t,m}$ be the covariate of interest that is $\sigma(X_{t-m+1}, \cdots, X_t)$ measurable. Then given $V_{t,m}$ we want to find the $\alpha$-lower quantile for the conditional random variable $T_{t+1,m}|V_{t,m}$ as our  VaR prediction for the future $m-$period returns. This is analogous to the one-sided prediction bounds described in Section \ref{sec:conjecture}, and we can use either QE or MFB to estimate the quantile of interest.

In our experiment, we consider   1-minute level return data during trading hours for our simulation. When extracting data, we remove those that are gathered from the first and last 5-minute in each trading day for quality control.
We use the past $m$-period realized volatility $V_{t,m} = \sum_{k=t-m+1}^{t} X_k^2$ as the covariate and study its impact on future $T_{t+1,m}$ via VaR. We take $T_{t+1,m}$ to be the worst $m$-period cumulative return for the future data $(X_{t+1},\cdots, X_{t+m})$, and our experiments show that by choosing appropriate $m$ there exhibits clear heteroscedascity in the paired data $(V_{t,m}, T_{t+1,m})$. Being able to capture the relationship between them will help individual traders to control their risk exposure, especially for day trading practitioners. 

We formally setup the experiment as in the following. Let $\{X_n\}_{n=1}^{N}$ be the observed data up to time $N$; our goal is to predict the VaR of $T_{N+1,m}$ conditioning on $V_{N,n}$. The regression model is trained on observed data from the past $\{(V_{i,m}, T_{i+1,m})\}_{i\in \mathcal I}$ where $\mathcal I\subset \{m,\cdots,N-m\}$ is some index set.
Because both $V_{i,m}$ and $T_{i+1,m}$ are aggregated statistics over $m-$periods, they each tends to have strong serial dependence for adjacent time indices which violates the i.i.d. assumption in our regression setup. 
To fix this issue, we only consider data from the time grid $\mathcal I = \{(2k + 1)m, k\in\mathbb N\}$, such that $(V_{(2k+1)m,m}, T_{(2k+1)m + 1,m})$ are based on non-overlapping blocks of data $\{(X_{2mk + 1}, \cdots, X_{2m(k+1)}), k\in\mathbb N \}$, thereby reducing the sample size 
by a factor of $2m$.

A new challenge arises from this procedure. 
On the one hand, as the dynamics of price process is time-varying,
to ensure the joint distribution of $(V_{i,m}, T_{i+1,m})$ does not change over time, the past data $\{X_n\}_{n=1}^N$ needs to be from a short time frame, e.g., within 2 weeks so that the data size is less than $N_0 = (6.5\times 60 - 10)\times 10 = 3800$.
On the other hand, through our observation, the heteroscedasticity phenomena are most pronounced when $m$ is larger than some threshold $m_0$, e.g. $m_0 = 30$, although the choice of $m_0$ varies for different stocks. By our data trimming procedure above, the sample size of the paired data is at most $|\mathcal I| = \lfloor \frac{N_0}{2m_0} \rfloor  = 63$.
This situation naturally falls into our limited data scenario introduced in Section \ref{sec: conditionality}. 

We apply QE, MFB as well as CP
using the nonparametric CDF estimator
to predict the future VaR 
on two stocks, AMC and GME listed on NYSE; these two stocks have recently been
very popular with day traders, and have exhibited 
abnormal stock price upsurges despite drops in revenue reported in their financial statements. Coupled with this phenomenon, the volatility and  fluctuations in the returns of these two stocks are higher than market average. It is therefore of practical value to devise a more formal testing procedure. 

We gathered data from April 10th to July 10th 2021 using the AlphaVantage stock APIs and performed post processing outlined in the above paragraphs to get pairs of $(V,T)$ data. Our VaR predictions are calculated for $\alpha \in \{0.01, 0.05, 0.1\}$ on $(V_{t,m},T_{t+1,m})$ with $m = 45$ for GME and $m=30$ for AMC, where the training data used for testing either have $|\mathcal I| = 34$ or $|\mathcal I| = 68$ prior to time $t$. 
We perform conjecture testing based on predicted VaRs  and future values $T_{t+1,m}$, with null hypothesis being $T_{t+1,m} > VaR_{pred}$.
This testing procedure is carried out consecutively with respect to $t$ so that we have a series of accept-reject results.
Lastly, we calculate the empirical acceptance rate based on the results and compare it with the nominal level, $1-\alpha$. We benchmark the methods by  closeness between the empirical acceptance rate and the nominal level.

The results of our experiments are presented in Table \ref{tab:simulation_GME} and \ref{tab:simulation_AMC}. 
The methods used for comparison are: na\"ive quantile estimation(QE); Model-free bootstrap with $L_1$ and $L_2$ predictors (MFB-$L_1$ and MFB-$L_2$ respectively); and (one-sided) distributional conformal prediction(CP). 

Based on the empirical results, both QE and MFB still exhibit undercoverage issue, while we find that overall MFB is able to boost the coverage of QE for both $|\mathcal I| = 34$ and $|\mathcal I| = 68$, exemplifying the result of Theorem \ref{thm:performance}; when sample size increases to $|\mathcal I| = 68$, the undercoverage issue is relieved to some level for both methods. 

Measured by closeness to the nominal level, for the AMC stock the performance of CP 
is superior to other participating methods, while for GME the MFB has the edge.
Notice that for the GME stock and some datapoints of the AMC stock, there are overcoverage issues for CP,  echoing our analysis for the behavior of CP under limited-data scenario in the previous numerical experiment.

\begin{table}[h]
	\centering
	\begin{tabular}{c c c  c  c}
		\toprule[0.25ex]
		\multicolumn{5}{c}{Ticker: GME}\\
		\hline 
		Sample size & Method &$\alpha = 0.01$ &$ \alpha = 0.05$ &$ \alpha = 0.1 $\\
		\hline
		\multirow{4}{*}{ $|\mathcal I| = 34$} & QE &95.9     & 94.1   &   88.2\\
		&MFB-$L^1$ &{97.2}    &94.1   &   {89.1}\\
		&MFB-$L^2$ &97.7     &94.1   &  {88.6} \\
		&CP & 99.5 & 97.7 & 91.4\\
		\midrule[0.2ex]
		\multirow{4}{*}{ $|\mathcal I| = 68$} & QE & {96.7}    & 94.1   &89.8   \\
		&MFB-$L^1$ &  97.3   & 95.7   &89.8  \\
		&MFB-$L^2$ & 98.6  &95.1    &89.8  \\
		&CP & 99.4 & 96.7 &92.5\\
		\bottomrule[0.25ex]
	\end{tabular}
	\caption{Empirical acceptance rate results for GME.  } 
	\label{tab:simulation_GME}
\end{table}

\begin{table}[h]
	\centering
	\begin{tabular}{c c c  c  c}
		\toprule[0.25ex]
		\multicolumn{5}{c}{Ticker: AMC}\\
		\hline 
		Sample size & Method &$\alpha = 0.01$ &$ \alpha = 0.05$ &$ \alpha = 0.1 $\\
		\hline
		\multirow{4}{*}{ $|\mathcal I| = 34$} & QE &95.1 & 90.8 & 86.5\\
		&MFB-$L^1$ & 95.1 & 91.6 & 87.6\\
		&MFB-$L^2$ & 95.4 & 91.4 &86.8 \\
		&CP &98.2 & 97.1 & 90.5\\
		\midrule[0.2ex]
		\multirow{4}{*}{ $|\mathcal I| = 68$} & QE & 95.8&91.7 & 87.9\\
		&MFB-$L^1$ &97.1    & 92.6    &   88.2\\
		&MFB-$L^2$ &  96.8  & 92.3   & 87.9  \\
		&CP &99.0 &95.8 & 89.8\\
		\bottomrule[0.25ex]
	\end{tabular}
	\caption{Empirical acceptance rate results for AMC. } 
	\label{tab:simulation_AMC}
\end{table}

\appendix

\section{Asymptotic analysis for CP}\label{appendix:A1}
Denote $\widehat U(x,y) = \widehat F(y|x)$, $\widehat U_i = \widehat U(X_i,Y_i)$, and $\widehat U_{(\alpha/2)}$, $\widehat U_{(1-\alpha/2)}$ the $\alpha/2$ sample quantiles of  $\{\widehat U_i\}$. Note that in \cite{distributionalconformal}, the conditional CDF is estimated with both $\{(X_i,Y_i)_{i=1}^n\}$ and $(X_{f},y)$, which results in an augmentation to $\widehat U_i$ that is denoted by $\widehat U_i'(X_f, y) = \widehat F_{n+1}(Y_i|X_i)$.

We assert that $\widehat U_i$ and $\widehat U_i'(X_f, y)$ are asymptotically equivalent. Specifically, 
let $ (X_{f} , Y_{f}) $ take on arbitrary values $ (x, y)$; then,  $\sup_{i} |\widehat U_i - \widehat U_i'(x, y)| = \mathcal{O}_p(1/n)$. As both $\widehat U_{(\alpha)}$ and $\widehat U_{(\alpha)}'$ converge to $\alpha$ in probability, $\widehat U_{(\alpha)} - \widehat U_{(\alpha)}' = o_P(1)$. Therefore, we can use $\widehat U_i$  instead for simplified asymptotic analysis.

Although in \cite{distributionalconformal}, the estimated ranks $\widehat U_i'$ are further transformed to $\widehat V_i$ with which the $p$-values are calculated, the key to generating prediction interval still relies on the ranks themselves. With this in mind,  the conditional coverage probability  approximately equals to

\begin{equation}\label{p_conformal}
	\mathbb P\left(\widehat U(X_{f}, Y_{f}) \in \left(\widehat U_{(\alpha/2)}, \widehat U_{(1-\alpha/2)}\right) |X_{f} = x_f\right).
\end{equation}
Let $\mathbb P_n$ and $\mathbb E_n$ denote the probability measure and expectation
respectively associated with data $(\mathbf X_n, \mathbf Y_n)$.
Then, by iterated expectation, (\ref{p_conformal}) equals to
\begin{equation*}
	\begin{split}	
		&\mathbb{E}_{y_f\sim F(\cdot|x_f)}\mathbb P_n\left(\widehat U(x_f,y_f) \in \left(\widehat U_{(\alpha/2)}, \widehat U_{(1-\alpha/2)}\right) |X_{f} = x_f\right)\\
		&= \mathbb{E}_{y_f\sim F(\cdot|x_f)}\mathbb E_n \mathbbm{1}\left\{\widehat U(x_f,y_f) \in \left(\widehat U_{(\alpha/2)}, \widehat U_{(1-\alpha/2)}\right) \right\}\\
		&=\mathbb E_n \mathbb{E}_{y_f\sim F(\cdot|x_f)} \mathbbm{1}\left\{y_f \in \left( \widehat F^{-1}_{y|x_f}(\widehat U_{(\alpha/2)}),   
		\widehat F^{-1}_{y|x_f}(\widehat U_{(1-\alpha/2)}) \right)\right\}\\
		&\approx \mathbb E_n F_{y|x_f}(\widehat F_{y|x_f}^{-1}(\widehat U_{(1-\alpha/2)})) - F_{y|x_f}(\widehat F_{y|x_f}^{-1}(\widehat U_{(\alpha/2)}))
	\end{split}
\end{equation*}
Similarly as before, denote $A_n' = \widehat F_{y|x_f}^{-1}(\widehat U_{(1-\alpha/2)}) - F_{y|x_f}^{-1}(\widehat U_{(1-\alpha/2)})$ and $B_n' = \widehat F_{y|x_f}^{-1}(\widehat U_{(1-\alpha/2)}) - F_{y|x_f}^{-1}(\widehat U_{(1-\alpha/2)})$ the errors of estimated quantiles; using Taylor expansion we get a similar expansion as in Equation (\ref{eq:coveragebias}), i.e., 
\begin{align}\label{eq:cvpconformal}
	\begin{split}
		&\mathbb E_n F_{y|x_f}(\widehat F^{-1}(\widehat U_{(1-\alpha/2)})) - F_{y|x_f}(\widehat F^{-1}(\widehat U_{(\alpha/2)})) \\
		&= (1-\alpha) + \mathbb E_n(f_{y|x_f}(F_{y|x_f}^{-1}(\widehat U_{(1-\alpha/2)})))\mathbb E_n(B_n') - 
		\mathbb E_n(f_{y|x_f}(F_{y|x_f}^{-1}(\widehat U_{(\alpha/2)})))\mathbb E_n(A_n')\\
		&\quad+\frac{1}{2} \mathbb E_n(f_{y|x_f}'(F_{y|x_f}^{-1}(\widehat U_{(1-\alpha/2)}))) Var_n(B_n') - \frac{1}{2}
		\mathbb E_n(f_{y|x_f}'(F_{y|x_f}^{-1}(\widehat U_{(\alpha/2)}))) Var_n(A_n') + \text{higher-order terms.}
	\end{split}	
\end{align}
Hence, the bias in coverage probability is
directly due to the bias and variance of the quantile estimator; in that sense,
QE and CP perform very similarly.

The most notable difference between QE and CP is the addition of a new data point $(X_f, y)$ in the original samples $\{(X_i, Y_i)\}_{i=1}^n$ to estimate $\widehat F$, which results in the augmentation $\widehat U_i'(X_f, y)$ from $\widehat U_i$. An intuitive interpretation of this is that the new procedure “overfits” on a future observation $(X_f, y)$ and may lead to better coverage with finite sample size. With exchangeability holding, the procedure has guaranteed unconditional validity, while conditional validity holds asymptotically. 

\subsection{Asymptotic analysis for MFB}\label{appendix:A2}

We adopt the usual notation, where $\mathbb P$ denotes the probability measure for the original probability space, and $\mathbb P^*$ the probability measure in the bootstrap world, i.e., conditioning on the data $(\mathbf X_n, \mathbf Y_n)$.
Given future covariate $x_f$, let $D_{x_f}$ denote the distribution of the predictive root $R_f = Y_f - \widehat Y_f$ under $\mathbb P$, and let $\widehat D^*_{x_f}$ and $\widehat D_{x_f}^{*-1}$
denote the distribution and quantile respectively of the bootstrap version $R_f^*$ under $\mathbb P^*$; also assume that both are continuous functions, i.e., that  $R_f$ is a
continuous random variable.  The bootstrap prediction interval 
$\mathcal{\widehat C}^{MF}(X_f)$ is $ (\widehat Y_f + \widehat D_{x_f}^{*-1}(\alpha/2), \widehat Y_f + \widehat D_{x_f}^{*-1}(1-\alpha/2))$, and 
is considered a random element under $\mathbb P$. The randomness of $\mathcal{\widehat C}^{MF}(X_f)$  
is due to variability across samples; given the observed sample, $\mathcal{\widehat C}^{MF}(X_f)$ 
is a fixed quantity. 

For the moment, let us focus our interest on the  coverage probability without conditioning on the
observed samples; hence the governing law is the original probability measure $\mathbb{P}$. In addition, we have the following lemma concerning the distribution of bootstrap samples.

\begin{lemma}\label{lemma:bootstrapdistribution}
	
	1. Under $\mathbb P^*$, the distribution of $X_i^*$ conditioning on $(\mathbf X_n, \mathbf Y_n)$ is $\frac{1}{n} \sum \mathbbm 1\{\cdot \leq X_i\}$; The distribution of $Y_i^*| X_i^* = X_k$ is $\widehat F_{y|X_k}$.\\
	2. Under $\mathbb{P}$, the distributions of $X_i^*$ and $Y_i^*$ are the expectations of their distributions under $\mathbb P^*$. i.e., 
	the distribution of $X_i^*$  is $F_{X}$, and $Y_i^*|X_i^* = x$ follows $\mathbb E_n \widehat F_{y|x}$. \\
	3. Under $\mathbb P^*$, $Y_f^*$ is independent of  $\widehat Y_f^*$, so
	the distribution of $R_f^* = Y_f^* - \widehat Y_f^*$ is the convolution between $Y_f^* \sim \widehat F_{y|x_f}$ and the distribution of $\widehat Y_f^* = g\left(\{X_i^*, Y_i^*\}_{i=1}^n\right)$.\\
	4. Under $\mathbb P$, $Y_f$ is independent of $\widehat Y_f$ as well, and the distribution $D_{x_f}$ is the convolution of $F_{y|x_f}$ and the distribution of $\widehat Y_f$.
\end{lemma}

Let $T_z f(x) = f(x + z)$ be the translation operator. For a CDF $F$, $(T_z F)^{-1}(u) = F^{-1}(u) - z $. 
Denote 
$\Delta_1^{(u)} = \widehat D_{x_f}^{*-1} (u) - \left(T_{ \mathcal J(\widehat F^*_{y|x_f}) }\widehat F_{y|x_f}\right)^{-1}(u)$; $\Delta_2 = \mathcal J(\widehat F^*_{y|x_f}) - \mathcal J(\widehat F_{y|x_f})$. Also, $A_n$ and $B_n$ are the same as in section \ref{sec: asymptotics}.

Note that  $\widehat Y_f$ and $R_f^*$ are dependent as they both depend  on the data, so the  previous route for coverage calculation can not be directly applied here. To carry out some asymptotic analysis, we  make the following extra assumptions:   $\widehat Y_f$ and $\widehat Y_f^*$ converge  (respectively) to $\mathcal J(F_{y|x_f})$ and $\mathcal J(\widehat F_{y|x_f})$, which are estimands of the two, 
in probability and in  $L^1$. Also assume $D_{x_f}$  and $\widehat D_{x_f}^{*}$ are uniformly continuous, so that by
Taylor expansion ($\delta-$method), the distributions $D_{x_f}(x)$ and $\widehat D_{x_f}^{*}(x)$ converge (respectively)  to $T_{\mathcal J(F_{y|x_f})}F_{y |x_f}(x)$ and $T_{ \mathcal J(\widehat F_{y|x_f}) }\widehat F_{y|x_f}(x)$, uniformly in $x$. As a consequence, $\forall u\in (0,1)$, $\widehat D^{*-1}_{x_f}(u) \overset{P^*}{\rightarrow} \left(T_{ \mathcal J(\widehat F_{y|x_f}) }\widehat F_{y|x_f}\right)^{-1}(u) = \widehat F_{y|x_f}^{-1}(u) - \mathcal J(\widehat F_{y|x_f})$, and thus $\Delta^{(u)}\overset{P^*}{\rightarrow} 0$.

Consider the following decomposition holds:
\begin{equation}\label{eq:error_decomposition}
	\begin{aligned}
		&\widehat Y_f + \widehat D^{*-1}_{x_f}(\alpha/2) = F^{-1}_{y|x_f}(\alpha/2)+ A_n + \Delta_1^{(\alpha/2)} + \Delta_2,\\
		&\widehat Y_f + \widehat D^{*-1}_{x_f}(1-\alpha/2) = F^{-1}_{y|x_f}(1-\alpha/2)+ B_n + \Delta_1^{(1-\alpha/2)} + \Delta_2\\
	\end{aligned}
\end{equation}
Denote the total error in (\ref{eq:error_decomposition}) $e_{total} = A_n + B_n + \Delta_1^{(\alpha/2)} + \Delta_1^{(1-\alpha/2)} + \Delta_2$. Now by iterated expectation and Taylor expansion, the coverage probability for the PI \eqref{eq:MFB RR PI} is
\begin{align}
	\begin{split}
		&\mathbb E_n\left(F_{y|x_f}(\widehat Y_f + \widehat D^{*-1}_{x_f}(1-\alpha/2)) - F_{y|x_f}(\widehat Y_f + \widehat D^{*-1}_{x_f}(\alpha/2))\right)\\
		&=\mathbb E_n F_{y|x_f}\left( F^{-1}_{y|x_f}(1-\alpha/2)+ A_n +  \Delta^{(1-\alpha/2)} + \Delta_2\right) - \mathbb E_n F_{y|x_f}\left( F^{-1}_{y|x_f}(\alpha/2)+ B_n +  \Delta_1^{(\alpha/2)} + \Delta_2 \right) \\
		& = (1-\alpha) +O(\mathbb E_n e_{total} + Var_n (e_{total})) + \text{higher-order terms}.
	\end{split}	
\end{align}

Despite different formalization, the above discussion shows the three methods have similar
asymptotic behavior for their respective  coverage probability when conditioning on the future covariate $x_f$ \color{black}, in the sense that the source of coverage bias is from the first two moments of the total estimation error in each method.

\section{Proof of Lemma \ref{lm:convolution}}

\begin{proof}
	First of all, we asssume (i) holds for all $u\in\mathbb R$. 
	By Taylor's Theorem,
	
	$$f(u + z) = f(u) + f'(u)z + \frac{f''(u)}{2} z^2 + R(z),$$
	where the remainder $R(z)$ has the  form $$R(z) = \frac{f^{(3)}(\xi(z))}{3!}z^3$$
	for some $\xi(z) \in (u, u + z).$
	Since $|f^{(3)}(\xi(z))|$ is bounded, we have $$\int_{\mathbb R} R(z) \phi_\sigma(z)\, dz = o(\sigma^2),$$
	where $\phi_\sigma (z)$ is the density function for the normal distribution $\mathcal N (0, \sigma^2)$. 
	
	Denote the density of $Y + \sigma Z$, where Y has density $f(u)$ and $\sigma Z \sim \mathcal N(0,\sigma^2)$ as $\tilde f$. Then by convexity, $\forall |u|> u_0$
	\begin{equation}\label{eq: inequal}
		\begin{split}
			\tilde f(u) = (f\star \phi_\sigma)(u) &= \int_{\mathbb{R}} f(u + z)\phi_\sigma(z) dz\\
			&=  \int_{\mathbb{R}} \left(f(u) + f'(u)z + \frac{f''(u)}{2} z^2 + R(z)\right)\phi_\sigma(z)\, dz\\
			&= f(u) + \frac{f^{''}(u)}{2}\sigma^2 + o(\sigma^2)\\
			& > f(u)
		\end{split}
	\end{equation}
	for small enough $\sigma^2$. As a consequence, $\forall u > u_0$, $\widebar{F\star\phi_{\sigma}}(u) = \displaystyle\int_{u}^{\infty} \tilde f(x)\, dx > \displaystyle\int_{u}^{\infty} f(x) \,dx = \widebar{F}(u)$. Analogously, $F\star\phi_{\sigma}(u) > F(u)$ for $u<-u_0$.

	In the more general case where (i) is satisfied on the restricted domain $|u| > u_0$, $\forall \epsilon>0$, we can choose $\sigma$ small enough where
	$\mathbb P(|\sigma Z| > u_0)<\epsilon$, so that $\forall |u| >2u_0$, $|u + \sigma Z|$ is concentrated within $(u_0,\infty)$ with  probability greater than $1-\epsilon$; therefore $f(u + \sigma Z)$ satisfies (i) and (ii) for $|u| > 2u_0$ with probability $1 - \epsilon$. By the same procedure of Equation (\ref{eq: inequal}), this gives 
	
	$$\tilde f(u) > f(u)(1-\epsilon) + \frac{f^{''}(u)}{2}\left(\sigma^2 - \int_{z>c_{\sigma, \epsilon/2}} z^2\phi_\sigma(z)\, dz\right) + o(\sigma^2)$$
	where $c_{\sigma, \epsilon/2}$ is the upper $\epsilon/2$ quantile for the $\mathcal N(0, \sigma^2)$ distribution.
	Since $\epsilon$ is arbitrarily small and $f$ is convex for $|u| > 2u_0$, we conclude that $\tilde f(u) > f(u)$ for  $|u| > 2u_0$. The remainder of the proof is the same as above. 
\end{proof}

Note that we may relax condition (i) in lemma \ref{lm:convolution} where $f^{''}$ exists and is absolutely continuous on $u > |u_0|$. Then $f^{(3)}$ exists and the remainder $R(z)$ has an integral form $R(z) = \displaystyle\int _{(u, u + z)}{\frac {f^{(3)}(t)}{2}}(u + z -t)^{2}\,dt,$ and with bounded $f^{(3)}$, $\displaystyle\int_{\mathbb R} R(z)\phi_\sigma(z)\, dz$ is also $o(\sigma^2)$.

\section{Proofs  of  Theorems \ref{thm:asymptotics},  \ref{thm:pertinence} and \ref{thm:performance}}\label{appendix:A3}
\label{sec:proofs2}

{\itshape {\textbf {Proof of Theorem \ref{thm:asymptotics}.} }} 
By Lemma 1.2.1 of \cite{Politis1999}, under (C1) with $x = x_f$, $\widehat F^{-1}_{y|x_f}(\alpha) \overset{P}{\rightarrow} F^{-1}_{y|x_f}(\alpha), \forall \alpha \in (0,1).$ As a consequence, the boundary points of  $\widehat{\mathcal C}_{1-\alpha}^{(QE)}(x_f)$ converges to those of $\mathcal C_{1-\alpha}(x_f) = \left(F^{-1}_{y|x_f}(\alpha/2), F^{-1}_{y|x_f}(1-\alpha/2)\right)$, then  we have that 
$\mathbb P_3\left(Y_f \in \widehat{\mathcal C}^{(QE)}_{1-\alpha}(X_f) \right) \overset{P}{\rightarrow} 1-\alpha.$

With further assumption that (C1) holds for all $x$, then $\widehat U(Y_i|X_i)$(and $\widehat U'(Y_i|X_i)$) has asymptotic uniform distribution in $[0,1]$. To see this, note that for $X_i = x$ fixed, $\widehat U(Y_i|X_i =x) - U(Y_i|X_i = x) \overset{P}{\rightarrow} 0,$ then 
$\mathbb P(\widehat U(Y_i|X_i =x) \leq u) \rightarrow u$, $\forall u\in[0,1]$. Now 
$$\mathbb P(\widehat U(Y_i|X_i)\leq u) = \mathbb E_{X_i}\mathbb E_{Y_i|X_i = x}\mathbb E_{n-1} I(\widehat U(Y_i|X_i = x)\leq u) = \mathbb E_{X_i}\mathbb P(\widehat U(Y_i|X_i =x) \leq u),$$
then by dominated convergence theorem, $\mathbb P(\widehat U(Y_i|X_i)\leq u) \rightarrow u$.
Then for CP, the analysis in section \ref{sec: asymptotics} shows that $\widehat{\mathcal C}_{1-\alpha}^{(CP)}(x_f)$ is asymptotically equivalent to $\widehat{\mathcal C}_{1-\alpha}^{(QE)}(x_f)$.

With the additional assumption $\widehat Y_f^* - \widehat Y_f \overset{P}{\rightarrow} 0$, $\widehat{\mathcal C}_{1-\alpha}^{(MF)}(x_f)$ is asymptotically equivalent to $\widehat{\mathcal C}_{1-\alpha}^{(QE)}(x_f)$, thus MFB is asymptotically valid. 
\vspace{0.2cm}

\noindent{\itshape {\textbf {Proof of Theorem \ref{thm:pertinence}.} }} 
We decompose $R_f$ and $R_f^*$ as follows.
\begin{eqnarray}
	R_f = Y_f - \widehat Y_f = \left(Y_f - \mathcal J(F_{y|x_f})\right) - \left(\widehat Y_f - \mathcal J(F_{y|x_f})\right) : = \epsilon_f + e_f;\\
	R_f^* = Y_f^* - \widehat Y_f^* = \left(Y_f^* - \mathcal J(\widehat F_{y|x_f})\right) - \left(\widehat Y_f^* - \mathcal J(\widehat F_{y|x_f})\right): = \epsilon_f^* + e_f^*.
\end{eqnarray}

with $\epsilon_f = Y_f - \mathcal J(F_{y|x_f})$, $\epsilon_f^* = Y_f^* - \mathcal J(\widehat F_{y|x_f})$; and $e_f =  \mathcal J(F_{y|x_f}) - \widehat Y_f$, $e_f^* =  \mathcal J(\widehat F_{y|x_f}) - \widehat Y_f^*$. 

Under assumption (C1), (2) in definition \ref{def:pertinence} is satisfied. Under assumption (C2), both $\tau_n e_f$ and $\tau_n e_f^*$ converge in distribution to normal distributions with equal variances, which shows (3) in definition \ref{def:pertinence}. 

Observe that $\epsilon_f$ only depends on $Y_f$, and the randomness of $e_f$ depends on the data  $(\mathbf X_n, \mathbf Y_n)\indep Y_f$, therefore $\epsilon_f\indep e_f$. Similarly, conditioning on the data, $\epsilon_f^* \indep e_f^*$ in the bootstrap world. This shows (4).

\vspace{0.2cm}

\noindent{\itshape {\textbf {Proof of Theorem \ref{thm:performance}.} }} 
Under scenario (A3), assumption (C2) and previous result, the distribution of the bootstrap root  ($R_f^*$) is asymptotically equivalent to the  convolution between $\widehat F_{y|x_f}$ and a Gaussian distribution $\mathcal  N(0, \tau_n^2)$. Then assumption (C3) guarantees that $\exists \alpha_0\in (0,1)$, $n_0\in \mathbb N$ such that $\forall \alpha < \alpha_0$ and $n >n_0$, $\widehat{\mathcal C}_{1-\alpha}^{(QE)}(x_f) \subseteq \widehat{\mathcal C}_{1-\alpha}^{(MF)}(x_f)$. This shows Equation (\ref{eq:MFBadvantage}) for $k = 3$. By Lemma \ref{lm:tower}, we can show (\ref{eq:MFBadvantage}) for the other two scenarios.

\section*{Acknowledgment}
Many thanks are due to E. Candes, V. Chernozhukov, L. Wasserman,  K. W\"uthrich, and Y. Zhang 
for discussions on conditionality, model-free methods and conformal prediction. 
Research partially supported by NSF grant DMS 19-14556 {and the Richard Libby graduate research award}.

\bibliography{reference}

\end{document}